\documentclass[preprint,12pt]{elsarticle}

\usepackage{amssymb}
\usepackage{amsmath}

\journal{Nuclear Physics B}
\usepackage[left=2.5cm, right=2.5cm, top=3cm, bottom=3cm]{geometry}
\begin{document}

\begin{frontmatter}



\title{Nonlinear Dynamics and Performance Optimization Based on Primary Resonance of an Electromechanically Coupled Magnetic Levitation Energy Harvester
}


\author{Jinhao~Xie~~~Chengkai~Yuan} 

\affiliation{organization={East China University of Science and Technology},
            addressline={No. 130, Meilong Road, Xuhui District}, 
            city={Shanghai},
            postcode={200237}, 
            country={China}}

\begin{highlights}

\item  A modified magnetic levitation system is developed by incorporating a more realistic energy harvesting circuit, enabling a better representation of practical operating conditions.

\item  An extended detuning formulation is introduced to capture the interaction between the harvesting circuit and the mechanical system, allowing the analysis of their coupled dynamic response under varying operating conditions.

\item  The proposed coupling mechanism effectively suppresses undesirable nonlinear behaviors such as chaos and multi-stability, resulting in a more stable and predictable energy harvesting process.

\item  A trade-off between energy harvesting efficiency and dynamical stability is identified, providing a new design perspective for practical energy harvesting systems.

\end{highlights}

\begin{abstract}
This research paper explores the potential of \textbf{nonlinear magnetic levitation systems} for energy harvesting by developing a modified system that incorporates a more realistic energy harvesting circuit, enabling a better representation of practical operating conditions.  
Methodologically, approximate solutions for the system dynamics were obtained using the \textbf{method of multiple scales}, complemented by numerical simulations to capture parameter variations visualized through phase planes and parameter variation plots. The results demonstrate that by adjusting capacitance to induce \textbf{internal and primary resonances}, an \textbf{extended detuning formulation} (utilizing parameters \(\sigma_3\) and \(\sigma_4\)) is innovatively introduced to capture the coupled dynamic interaction between the harvesting circuit and the mechanical system under diverse conditions. Periodic variations in circuit charge and intermediate magnet displacement were thoroughly analyzed. Comparisons with legacy models demonstrate that the proposed \textbf{coupling mechanism} effectively suppresses undesirable nonlinear behaviors, such as chaos and multi-stability, resulting in a more stable and predictable energy harvesting process. Ultimately, \textbf{a critical trade-off between energy harvesting efficiency and dynamical stability} is identified, providing a valuable new design perspective for practical energy harvesting systems.
\end{abstract}


\end{frontmatter}



\section{Introduction}
\label{sec1}
In recent years, technologies such as wireless sensor networks and microelectromechanical systems (MEMS) have achieved tremendous progress and development while simultaneously placing higher demands on energy supply \cite{3}. Traditional battery powering requires periodic replacement due to limited lifespan and may also cause environmental pollution. Consequently, for devices where battery replacement is impractical, attention has shifted toward seeking long-term power sources. Methods that harvest minute amounts of energy from the environment to charge batteries or directly power devices like sensors, thereby extending operational duration and enabling long-term, reliable operation, have garnered widespread attention.

This has positioned energy harvesting as a highly promising field within sustainable power generation. Its objective is to utilize renewable energy sources within the environment to meet the world's growing demand for clean, efficient electricity and promote sustainable development. Renewable energy refers to resources that can be replenished or regenerated from natural sources, including hydropower, geothermal energy, solar energy, wind energy, biomass, and tidal and wave energy \cite{8,9}. Among these, one of the most promising sources for energy recovery is the periodic vibrations generated by rotating machinery or engines \cite{10}.

Among numerous energy harvesting systems, magnetic levitation systems have garnered significant attention for their ability to convert mechanical vibrations into electrical energy\cite{1,2,333,444,555}. However, practical applications of these systems often deviate from ideal conditions, giving rise to complex issues that remain partially unresolved.

The reference\cite{11} explores the potential for energy extraction from the nonlinear vibrations of magnetic levitation. This work developed a novel device featuring tunable resonance and investigated the Duffing oscillator characteristics of the system under static and dynamic loading. The study posited that leveraging nonlinear phenomena could enhance the efficiency of energy harvesting devices. Reference\cite{12} examined magnetic levitation technology as a vehicle propulsion method, laying the groundwork for understanding its potential applications in high-speed train transportation. In 2017, Reference\cite{13} conducted comparative simulation and experimental studies on vibration energy harvesting using electromagnetic and piezoelectric seismic sensors. References\cite{14,15} further investigated the feasibility of converting vibrations into electrical energy using magnetic levitation collectors. Employing numerical simulations and experimental measurements, these studies emphasized the importance of adjusting electrical parameters to achieve optimal efficiency. The interaction between mechanical and electrical components was characterized through electromechanical coupling, with comparisons made between fixed-value, linear, and nonlinear coupling models. Results indicated that nonlinear electromechanical coupling models were better suited to high-vibration scenarios, where adjusting magnetic coil positions could control nonlinear resonance and enable energy recovery.

Reference\cite{16} establishes a mathematical model for a novel magnetic levitation device, accounting for non-ideal excitation of the main system by an electric vibrator. The magnetic levitation system here adopts the Duffing oscillator form, necessitating nonlinear analysis to investigate the energy harvesting potential of this nonlinear system. A multiscale approach is applied to study the vibration modes of the coupled system, highlighting its non-ideality and nonlinear phenomena. Results obtained using a fixed-step fourth-order Runge-Kutta numerical method are presented via phase plane plots, Bode plots, and parameter variations. These findings reveal multiple behaviors in magnet motion (periodic, quasi-periodic, and chaotic), coexisting attractors with high sensitivity to initial conditions, and intriguing outcomes for maximum average power: high and continuous energy acquisition within periodic and chaotic regions.

Subsequent reference\cite{17} further investigates the model from reference\cite{16}, emphasizing the impact of internal resonance between the excitation circuit and magnet motion on energy harvesting efficiency. It additionally discusses chaotic phenomena near excitation voltage amplitudes of E=3. The simulation section achieves internal resonance by varying capacitance, visually demonstrating magnet motion and power harvesting under resonance conditions. Shortly thereafter, the same research team published reference\cite{18}, focusing on the energy harvesting potential of this system when accounting for the weight of an intermediate suspended mass. This was achieved by introducing a new factor—the weight of the intermediate mass—through an axial translation geometric model. The models in references\cite{16}, \cite{17}, and \cite{18} all emphasize linking the force equations governing the central magnet's motion with the excitation circuit, highlighting the coupling between the excitation circuit and the magnet.

As early as 2019, reference\cite{19} investigated a similar electromagnetic energy harvester model. It employed the harmonic balance method to evaluate dynamic response and energy harvesting performance. Parametric analysis assessed the effects of damping, electromagnetic parameters, and coupling coefficients. The theoretical section emphasized linking the force equations governing the central magnet's motion with the energy harvesting circuit, highlighting the coupling between the circuit and the magnet.

This paper primarily integrates models from references \cite{16,17,18,19} to propose a more complex yet realistic and comprehensive nonlinear magnetic levitation model. \textbf{This work simultaneously considers the coupling between the excitation circuit, the energy harvesting circuit, and the mechanical system}, with the primary innovation being the correction of the crude modeling of the energy harvesting mechanism in the reference's\cite{16} model, which may fails to reflect specific application scenarios. The reference's\cite{16} model does not account for the counteraction of current variations in the energy harvesting coil on magnet motion and heavily relies on internal resonance to enhance energy harvesting efficiency. \textbf{This paper examines the influence of RLC energy harvesting circuits, with a stronger emphasis on enhancing the magnitude of the harvesting circuit current through primary resonance, serving as a crucial supplement to energy harvesting scenarios utilizing internal resonance}. Building upon this, the paper also discusses the impact of varying damping, excitation voltage, and intermediate magnet mass on energy harvesting efficiency, alongside stability analysis of magnet motion under the new model. \textbf{Furthermore, considering the coupling between the collection circuit and the mechanical system plays a positive role in increasing the stability of the system's energy harvesting and reducing the occurrence of chaotic phenomena}.

The structure of this paper is as follows: The preceding introduction constitutes the first part. The second part details the model of the system under study—the circuitry for energy harvesting. This section also introduces the excitation circuit and discusses the derivation process of the modeling equations. \textbf{Subsequently, the third part provides a detailed analysis of the multiple scales method for the normalized system equations,deriving the approximate response equations for the system}.The fourth part presents charts of numerical and analytical results, incorporating experimental data from references \cite{16,17,18,19}. It elaborates on how varying mechanical-electrical parameters—\textbf{such as mass, excitation voltage, and inductance—influence the excitation circuit current and the intermediate magnet's amplitude}, demonstrating the energy harvesting potential of the innovative system proposed. Furthermore, by comparing the model with data from prior references, the improvements in energy harvesting behavior achieved by the new model are highlighted, \textbf{and analysis is conducted using dynamical system methods such as phase portraits, Poincare maps, and bifurcation diagrams.}

\begin{figure}[htbp]
\centering
\centerline{\includegraphics[width=0.8\linewidth]{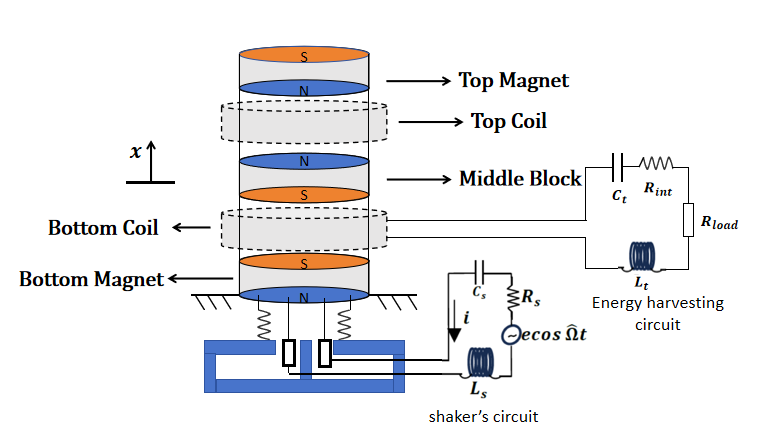}}
\caption{Experimental apparatus}
\label{fig1}
\end{figure}

\section{magnetic levitation system}

Figure \ref{fig1} illustrates the magnetic levitation structure connected to the energy harvesting device. In this design, the upper, middle, and lower magnets form the core assembly. When the excitation circuit is activated, the central magnet generates oscillations that alter the magnetic field, inducing an electric current in the harvesting coil. The magnetic levitation Faraday base is realized through an electromagnetic exciter controlled by an excitation circuit within an RLC framework, where $L_s$, $R_s$, and $C_s$ correspond to the inductance, resistance, and capacitance of the excitation circuit, respectively. To further investigate magnetic levitation energy harvesting behavior, this paper addresses potential capacitive and inductive components within the harvesting circuit. The purely resistive harvesting circuit governed by Kirchhoff's laws is \textbf{modified into an RLC-based harvesting circuit to better align with practical applications}. Here, $L_t$, $R_t$, and $C_t$ represent the inductance, resistance, and capacitance of the harvesting circuit, respectively.

\subsection{Original System Equations}
\label{subsec1}

The system equations of motion derived from the Euler-Lagrange equations are detailed in Ref.\cite{16}. The dynamic expressions proposed in this reference can be expressed as:
\begin{equation}
    m\ddot{x} + mg + k_1x + k_3x^3 + S_1\dot{q}_1 + C_{me}\dot{x} = 0
\end{equation}
\begin{equation}
   L_s\ddot{q}_1 - S_1\dot{x} + \frac{q_1}{C_s} + R_s\dot{q}_1 = e\cos(\widehat{\Omega}t)
\end{equation}
Here, $m$ denotes the mass of the central magnet, $x$ represents the displacement of the central magnet, $q_1$ denotes the charge in the excitation circuit, $k_1$ and $k_3$ denote the magnetic force parameters, $S_1$ denotes the coupling coefficient between the mechanical system and the excitation circuit, $C_{me}$ represents the sum of the mechanical damping coefficient and the electrical damping coefficient, $e$ indicates the amplitude of the AC excitation voltage, $\widehat{\Omega}$ denotes the angular frequency of the power source, and $g$ represents the gravitational acceleration.

The power expression for the energy harvesting circuit of this system is:
\begin{equation}
   P = i^2 R_{\text{load}} = \left( \frac{c_e}{\alpha} \dot{x} \right)^2 R_{\text{load}} 
\end{equation}
Here, $\alpha$ denotes the electromechanical interaction parameter, $i$ represents the current, and $R_{\text{load}}$ is the external resistance of the energy harvesting circuit.

\subsection{New System Equations}

This paper extends the system equations of the old system by adding analytical equations for the energy harvesting circuit. The generalized kinetic energy is derived from the motion of the intermediate magnet and the current flow in both circuits. The calculation of the mutual magnetic force between magnets employs the same method as in reference\cite{16}, detailed in Appendix A. Notably, after considering the modification of the harvesting circuit to an RLC circuit, \textbf{this paper simultaneously incorporates electromagnetic induction between the harvesting circuit and the mechanical system, establishing corresponding coupling coefficients to better align the energy harvesting behavior of the magnetic levitation system with practical conditions}. This refinement addresses ambiguities in the original magnetic levitation system's power harvesting calculations, enhancing its research value. Accordingly, the system equations under the new model are established as follows:

\begin{figure}[t]
\centering
\centerline{\includegraphics[width=1.2\linewidth, height=4cm]{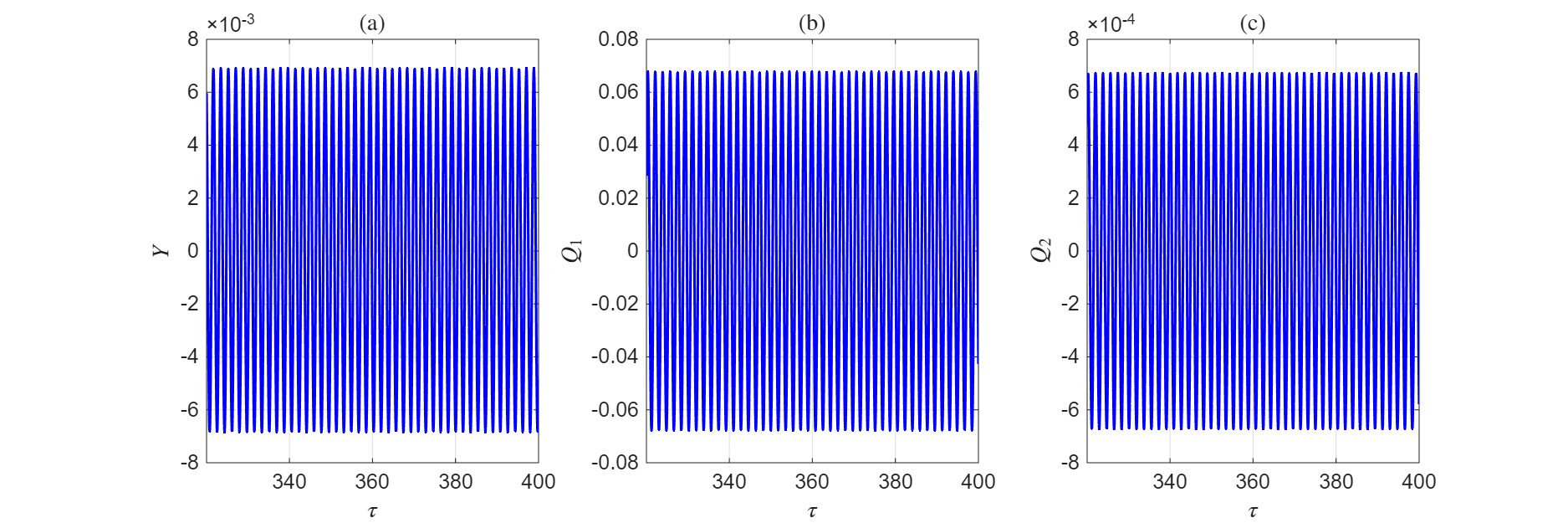}}
\caption{Response of three primary parameters to normalized time $\tau$ under internal resonance conditions}
\label{fig2}
\end{figure}

\begin{equation}
    m\ddot{x} + mg + k_1x + k_3x^3 + S_1\dot{q}_1 + S_2\dot{q}_2 + C_{me}\dot{x} = 0 \label{M}
\end{equation}
\begin{equation}
   L_s\ddot{q}_1 - S_1\dot{x} + \frac{q_1}{C_s} + R_s\dot{q}_1 = e\cos(\widehat{\Omega}t) \label{N}
\end{equation}
\begin{equation}
   L_t\ddot{q}_2 - S_2\dot{x}+\frac{q_2}{C_t} + R_t\dot{q}_2 = 0\label{A}
\end{equation}
Here, $m$ denotes the mass of the central magnet, $x$ represents the displacement of the central magnet, $q_1$ denotes the charge in the excitation circuit, $q_2$ denotes the charge in the collection circuit, and $k_1$ and $k_3$ denote the magnetic force parameters. $S_1$ and $S_2$ denote the coupling coefficients between the mechanical system and the excitation circuit, and between the mechanical system and the harvesting circuit, respectively, characterizing the intensity of energy exchange between the mechanical and electrical systems. $C_{me}$ denotes the sum of the mechanical damping coefficient and the electrical damping coefficient (where the mechanical damping component primarily arises from unavoidable factors such as mechanical friction and air resistance, while the electrical damping mainly stems from Joule heating generated by resistors in the energy harvesting circuit, characterizing the ability to transfer mechanical energy into electrical energy), $e$ represents the amplitude of the AC excitation voltage, $\widehat{\Omega}$ denotes the circular frequency of the power source, and $g$ represents the gravitational acceleration. Additionally, the displacement of the magnet's equilibrium position due to its weight must be considered. $Y_0$ denotes the equilibrium position when the magnet is at rest, expressed as:
\begin{equation}
   mg = k_1 Y_0 + k_3 Y_0^3 \label{0}
\end{equation}

\begin{figure}[t]
\centering
\centerline{\includegraphics[width=1.2\linewidth, height=4cm]{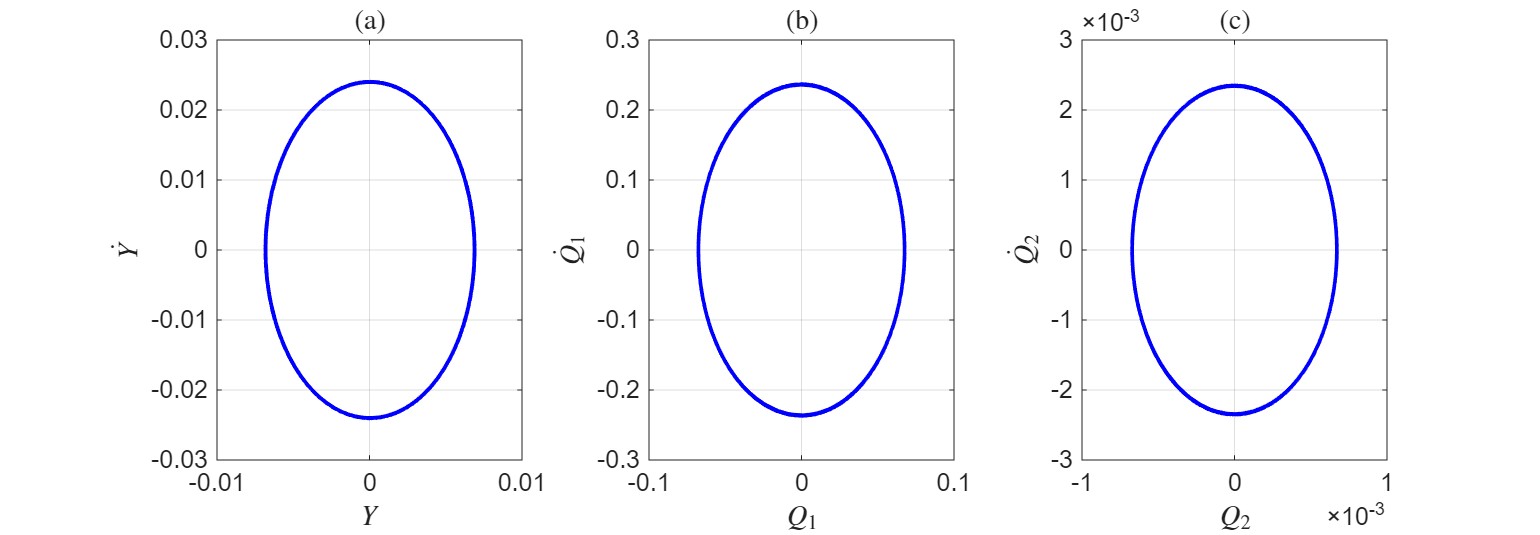}}
\caption{Stable phase plane of three primary parameters under internal resonance conditions}
\label{fig3}
\end{figure}

Set $\omega_0^2 = \frac{3k_3 Y_0^2 + k_1}{m}, \tau = \omega_0 t, Y = \frac{x}{x_0}, Q_1 = \frac{q_1}{q_0}, Q_2 = \frac{q_2}{q_0}$,\textbf{To normalize the above system of equations and eliminate the dimensional parameters for analysis}:

\begin{equation}
 \ddot{Y} + Y - W_2Y^2 + W_3Y^3 + \alpha_1\dot{Q}_1 + \alpha_2\dot{Q}_2 + \alpha_3\dot{Y} = 0 \label{1}
\end{equation}
\begin{equation}
  \ddot{Q}_1 - \alpha_4\dot{Y} + W_4Q_1 + \alpha_5\dot{Q}_1 = E\cos({\Omega}t)\label{2}
\end{equation}
\begin{equation}
   \ddot{Q}_2 + \beta_1\dot{Q}_2+ \beta_3{Q}_2- \beta_2\dot{Y} = 0\label{3}
\end{equation}
The parameters are defined as follows:
\begin{equation*}
\begin{aligned}
&W_2= \frac{3k_3 Y_0 x_0}{m \omega_0^2}, \quad W_3 = \frac{k_3 x_0^2}{m \omega_0^2}, \quad W_4 = \frac{1}{\omega_0^2 L_s C_s}, \quad \alpha_1 = \frac{S_1 q_0}{m \omega_0 x_0}, \\& \alpha_2 = \frac{S_2 q_0}{m \omega_0 x_0}, \quad \Omega = \frac{\hat{\Omega}}{\omega_0}, \quad \alpha_3 = \frac{C_{me}}{m \omega_0}, 
\quad \alpha_4 = \frac{S_1 x_0}{L_s \omega_0 q_0}, \quad \alpha_5 = \frac{R_s}{L_s \omega_0}, \\& E = \frac{e}{L_s \omega_0^2 q_0}, \quad \beta_1 = \frac{R_t}{L_t \omega_0}, \quad \beta_2 = \frac{S_2 x_0}{L_t \omega_0 q_0}, \quad \beta_3 = \frac{1}{\omega_0^2 L_t C_t}
\end{aligned}
\end{equation*}

The energy harvesting power under the new system is:
\begin{equation}
P =R_{\text{load}}i^2= R_{\text{load}}(\dot{Q}_2)^2
\end{equation}
Here, $R_{\text{load}}$ represents the load resistance of the energy harvesting circuit. Next, we will employ a perturbation method known as the multiscale approach to derive an approximate solution for this system.

\begin{figure}[t]
\centering
\centerline{\includegraphics[width=1.2\linewidth, height=5cm]{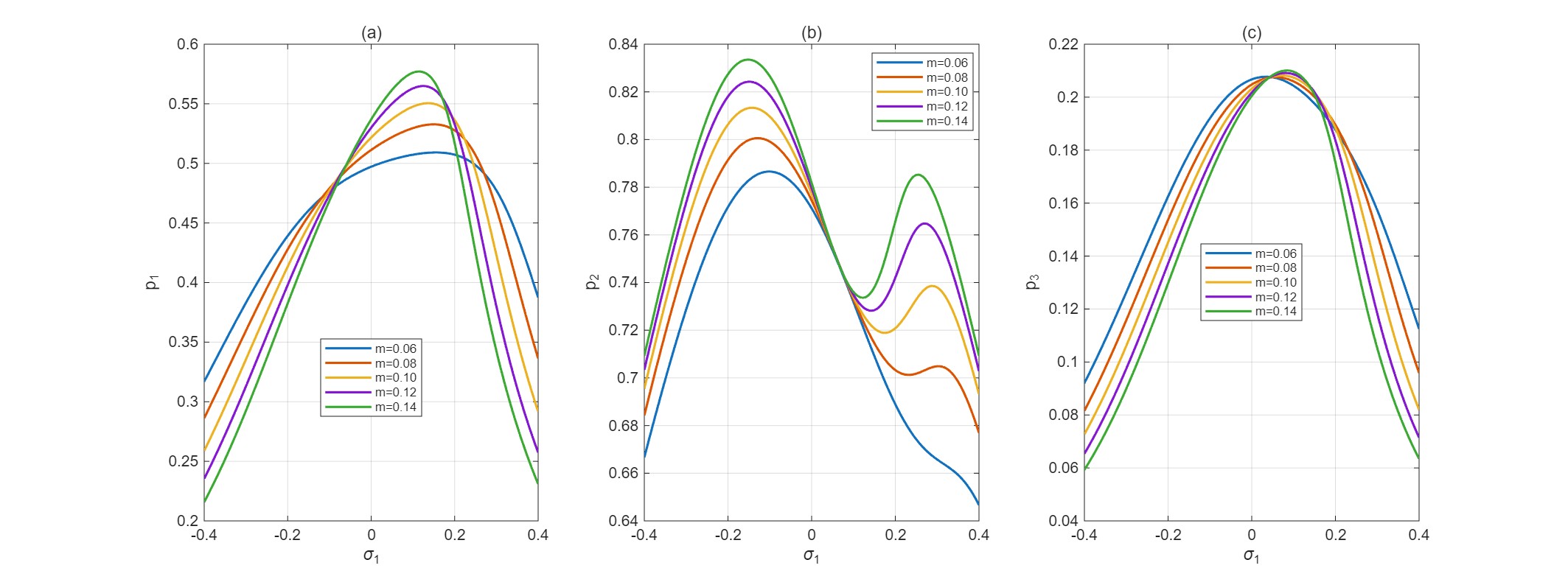}}
\caption{Amplitude of Three Primary Parameters at Different m Values Under Internal Resonance Conditions as a Function of $\sigma_1$}
\label{fig4}
\end{figure}

\section{Solve the approximate solution by using the multiple scales method }
\subsection{Internal Resonance Conditions}

The following section will employ a multiscale analysis approach to solve the differential equations. Here, a set of scaling parameters is introduced so that all parameters ultimately appear together in the perturbation equations.
\begin{equation*}
W_2 = \varepsilon \hat{W_2},  W_3 = \varepsilon \hat{W_3}, \alpha_i = \varepsilon\hat{\alpha_i} \ (i=1,2,3,4,5),E = \varepsilon \hat{E}, \beta_i=\varepsilon\hat{\beta_i} \ (i=1,2,3)
\end{equation*}
Let $T_0=\tau$ and $T_1=\varepsilon{\tau}$. Expand the first- and second-order derivative operators as power series in $\varepsilon$:
\begin{equation}
\begin{aligned}
\frac{d}{d\tau} &= D_0 + \varepsilon D_1 + O(\varepsilon^2)
\\ \frac{d^2}{d\tau^2} &= D_0^2 + 2\varepsilon D_0 D_1 + O(\varepsilon^2) \label{4}
\end{aligned}
\end{equation}
Expand $Y, Q_1, Q_2$ as power series in $\varepsilon$:
\begin{equation}
\begin{aligned}
Y &= u_0(T_0, T_1) + \varepsilon v_0(T_0, T_1) + O(\varepsilon^2)
\\
Q_1 &= u_1(T_0, T_1) + \varepsilon v_1(T_0, T_1) + O(\varepsilon^2)
\\
Q_2 &= u_2(T_0, T_1) + \varepsilon v_2(T_0, T_1) + O(\varepsilon^2) \label{5}
\end{aligned}
\end{equation}
\begin{figure}[t]
\centering
\centerline{\includegraphics[width=1.2\linewidth, height=4cm]{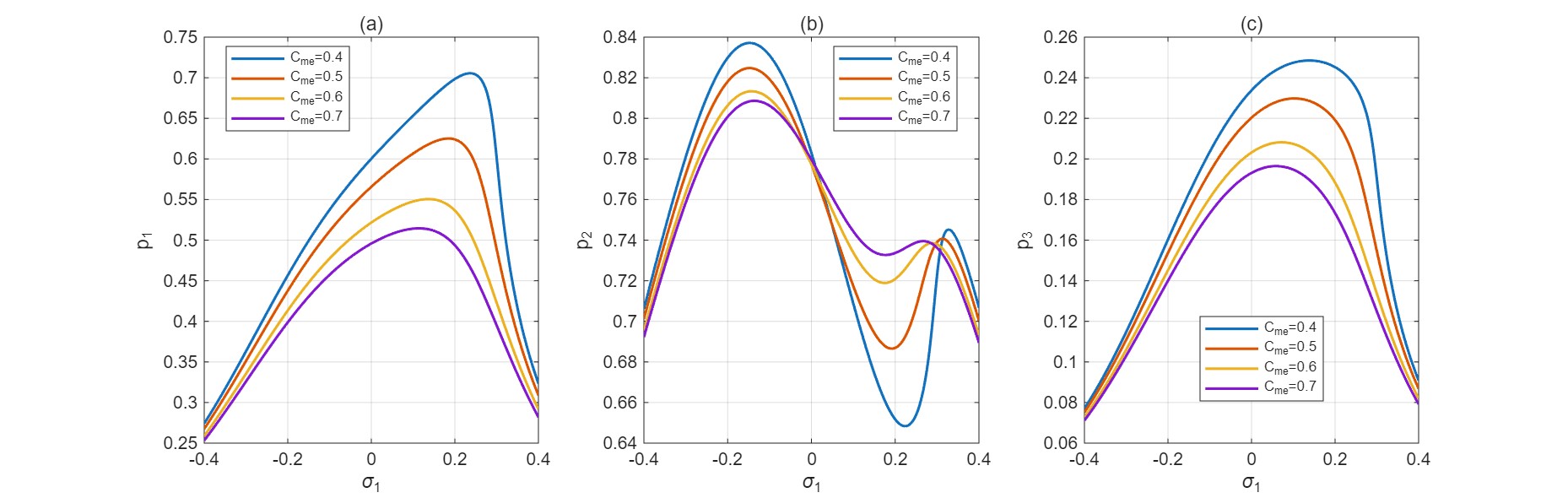}}
\caption{Amplitude of Three Primary Parameters at Different $C_{me}$ Values Under Internal Resonance Conditions as a Function of $\sigma_1$}
\label{fig5}
\end{figure}
Substituting equations(\ref{4}-\ref{5}) into equations(\ref{1}-\ref{3}) yields:\\
$O(\varepsilon^0)$:
\begin{equation}
\begin{aligned}
&D_0^2 u_0 + u_0 = 0
\\
&D_0^2 u_1 + W_4 u_1 = 0
\\
&D_0^2 u_2 + \beta_3 u_2= 0 \label{6}
\end{aligned}
\end{equation}

The final solutions are:
\begin{equation}
\begin{aligned}
&u_0 = A_1(T_1) e^{i T_0} + c.c.
\\
&u_1 = A_2(T_1) e^{i\sqrt{W_4} T_0} + c.c.
\\
&u_2 = {A_3(T_1)} e^{i\sqrt{\beta_3} T_0}+ c.c.\label{7}
\end{aligned}
\end{equation}
where $c.c.$ denotes the conjugate complex number.\\
$O(\varepsilon^1)$:
\begin{equation}
\begin{aligned}
&D_0^2 v_0 + 2 D_0 D_1 u_0 + v_0 - W_2 u_0^2 + W_3 u_0^3 + \alpha_1 D_0 u_1 + \alpha_2 D_0 u_2 + \alpha_3 D_0 u_0 = 0\\
&D_0^2 v_1 + 2 D_0 D_1 u_1 - \alpha_4 D_0 u_0 + W_4 v_1 + \alpha_5 D_0 u_1 = E \cos({\Omega} T_0)\\
&D_0^2 v_2 + 2 D_0 D_1 u_2+ \beta_1 D_0u_2 + \beta_3 v_2- \beta_2 D_0u_0= 0\label{8}
\end{aligned}
\end{equation}
We represent $A_i \ (i=1,2,3)$ using exponential complex numbers:
\begin{equation}
A_i = \frac{1}{2} p_i e^{i q_i} \label{8.5}
\end{equation}
\begin{figure}[t]
\centering
\centerline{\includegraphics[width=1.2\linewidth, height=4cm]{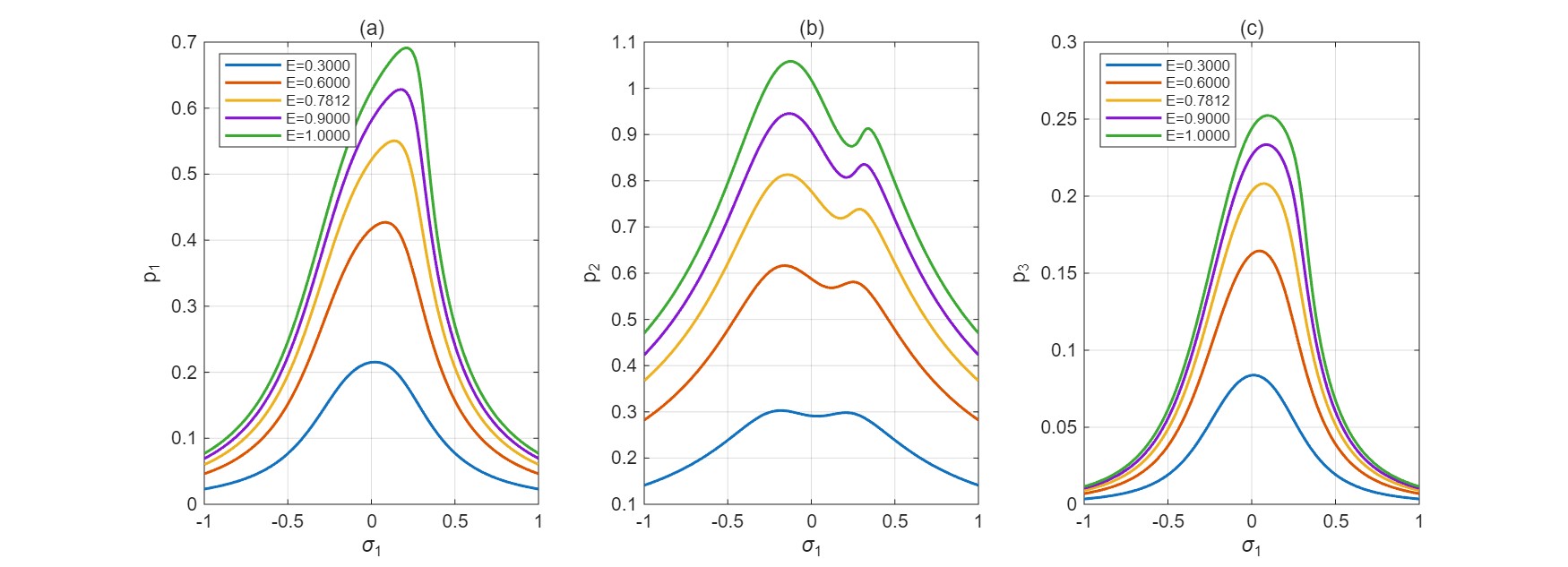}}
\caption{Amplitude of Three Primary Parameters at Different $E$ Values Under Internal Resonance Conditions as a Function of $\sigma_1$}
\label{fig6}
\end{figure}
To investigate the energy harvesting behavior of the system under different conditions, three parameters $\sigma_1$, $\sigma_2$, and $\sigma_3$ are introduced here:
\begin{equation}
\begin{aligned}
{\Omega} = \sqrt{W_4} + \varepsilon \sigma_1,  
\sqrt{W_4} = 1+ \varepsilon \sigma_2,  
\sqrt{\beta_3} =1+ \varepsilon \sigma_3 \label{9}
\end{aligned}
\end{equation}
where $\sigma_1$ is the primary resonance parameter, and $\sigma_2$ and $\sigma_3$ are the internal resonance parameters. Parameter $\sigma_1$ represents the closeness of $\sqrt{W_4}$ to $\Omega$, while parameters $\sigma_2$ and $\sigma_3$ represent the closeness of $\sqrt{W_4}$ and $\sqrt{\beta_3}$ to 1.

Substituting equations (\ref{7}, \ref{9}) into equation (\ref{8}) and rearranging the equation yields:
\begin{equation}
\begin{aligned}
(D_0^2 +1) v_0&=\left(-2iA_1' - i\alpha_1 \sqrt{W_4} A_2 e^{i\sigma_2 T_1} - i\alpha_2 \sqrt{\beta_3} A_3 e^{i\sigma_3 T_1} - i\alpha_3 A_1 - 3W_3 A_1^2 \overline{A_1}\right)e^{iT_0} + c.c.\\
(D_0^2 +W_4)v_1&=\left( i\alpha_4 A_1 e^{-i\sigma_2 T_1} - 2i\sqrt{W_4} A_2' - i\alpha_5 \sqrt{W_4} A_2 + \frac{1}{2} E e^{i\sigma_1 T_1} \right) e^{i\sqrt{W_4} T_0} + c.c.\\
(D_0^2 +\beta_3)v_2&=\left(-2i\sqrt{\beta_3} A_3' - i\beta_1 \sqrt{\beta_3} A_3 + i\beta_2 A_1 e^{-i\sigma_3 T_1}\right)e^{i\beta_3 T_0} + c.c. \label{10}
\end{aligned}
\end{equation}
\begin{figure}[t]
\centering
\centerline{\includegraphics[width=1.2\linewidth, height=4cm]{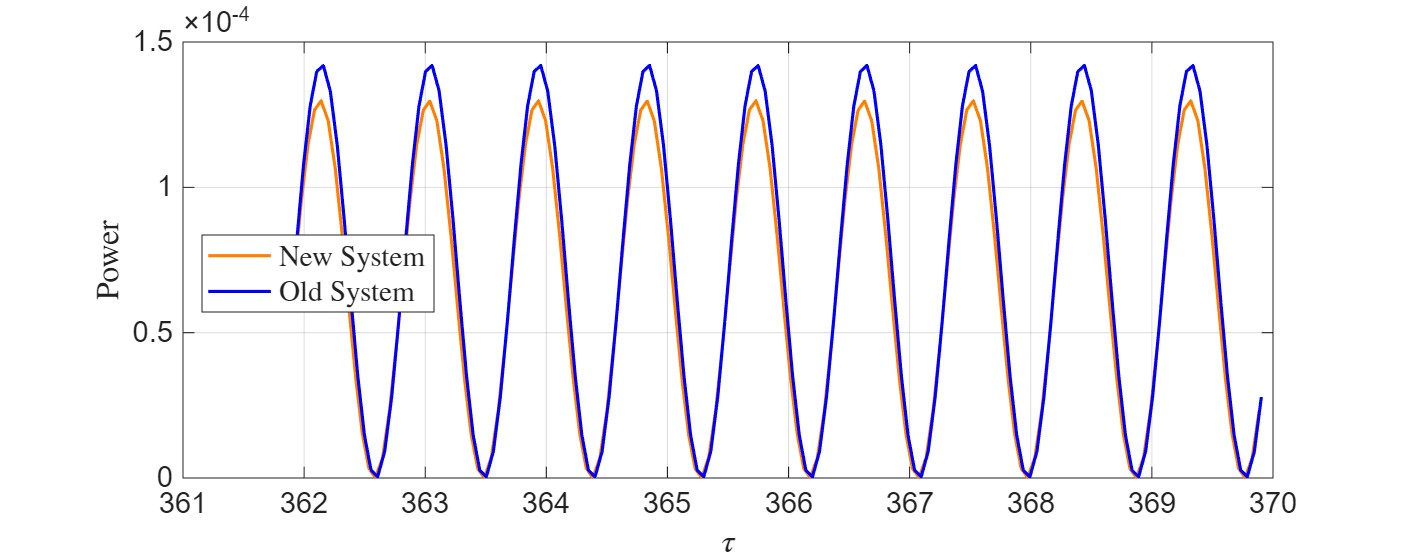}}
\caption{Power collected during internal resonance in the old system equation and primary resonance in the new system equation as a function of normalized time $\tau$}
\label{fig文献功率}
\end{figure}
Note that the right-hand sides of the three equations in (\ref{10}) must not contain terms with the same circular frequency as the corresponding solution to equation (\ref{7}), otherwise unbounded solutions will be obtained. This clearly does not match the actual response solution of the system. Therefore, to obtain stable solutions $v_0, v_1, v_2$, the part inside the brackets in equation (\ref{10}) must be zero:
\begin{equation}
\begin{aligned}
-&2iA_1' - i\alpha_1 \sqrt{W_4} A_2 e^{i\sigma_2 T_1} - i\alpha_2 A_3 \sqrt{\beta_3} e^{i\sigma_3 T_1} - i\alpha_3 A_1 - 3W_3 A_1^2 \overline{A_1} = 0\\
&i\alpha_4 A_1 e^{-i\sigma_2 T_1} - 2i\sqrt{W_4} A_2' - i\alpha_5 \sqrt{W_4} A_2 + \frac{1}{2}E e^{i\sigma_1 T_1} = 0\\
-&2i\sqrt{\beta_3} A_3' - i\beta_1 \sqrt{\beta_3} A_3 + i\beta_2 A_1 e^{-i\sigma_3 T_1}=0\\ \label{11}
\end{aligned}
\end{equation}
Three parameters are defined for subsequent derivation:
$\gamma_1 = q_2 + \sigma_2 T_1 - q_1$,
$\gamma_2 = \sigma_1 T_1 - q_2$,
$\gamma_3 = \sigma_3 T_1 + q_3 - q_1$.
Substituting equations (\ref{8.5}) and (\ref{9}) into equation (\ref{11}) yields:
\begin{equation}
\begin{aligned}
2p_1' &= -\alpha_1 \sqrt{W_4} p_2 \cos\gamma_1 - \alpha_2 \sqrt{\beta_3} p_3 \cos\gamma_3 - \alpha_3 p_1 
\\
2p_1 (\gamma_1' + \gamma_2') &= 2p_1 (\sigma_1 + \sigma_2) + \alpha_1 \sqrt{W_4} p_2 \sin\gamma_1 + \alpha_2 \sqrt{\beta_3} p_3 \sin\gamma_3 - \frac{3}{4} W_3 p_1^2 \overline{p_1} 
\\
2\sqrt{W_4} p_2' &= \alpha_4 p_1 \cos\gamma_1 - \alpha_5 \sqrt{W_4} p_2 + E \sin\gamma_2
\\
2p_2 \sqrt{W_4} \gamma_2' &= \alpha_4 p_1 \sin\gamma_1 + 2\sqrt{W_4} p_2 \sigma_1 + E \cos\gamma_2 
\\
2\sqrt{\beta_3} p_3' &= -\beta_1 \sqrt{\beta_3} p_3 + \beta_2 p_1 \cos\gamma_3 
\\
2p_3 \sqrt{\beta_3} (\gamma_3' - \gamma_1' - \gamma_2') &= -2p_3 \sqrt{\beta_3} (\sigma_1 + \sigma_2 - \sigma_3) - \beta_2 p_1 \sin\gamma_3\label{12}
\end{aligned}
\end{equation}
\begin{figure}[t]
\centering
\centerline{\includegraphics[width=1.0\linewidth, height=4cm]{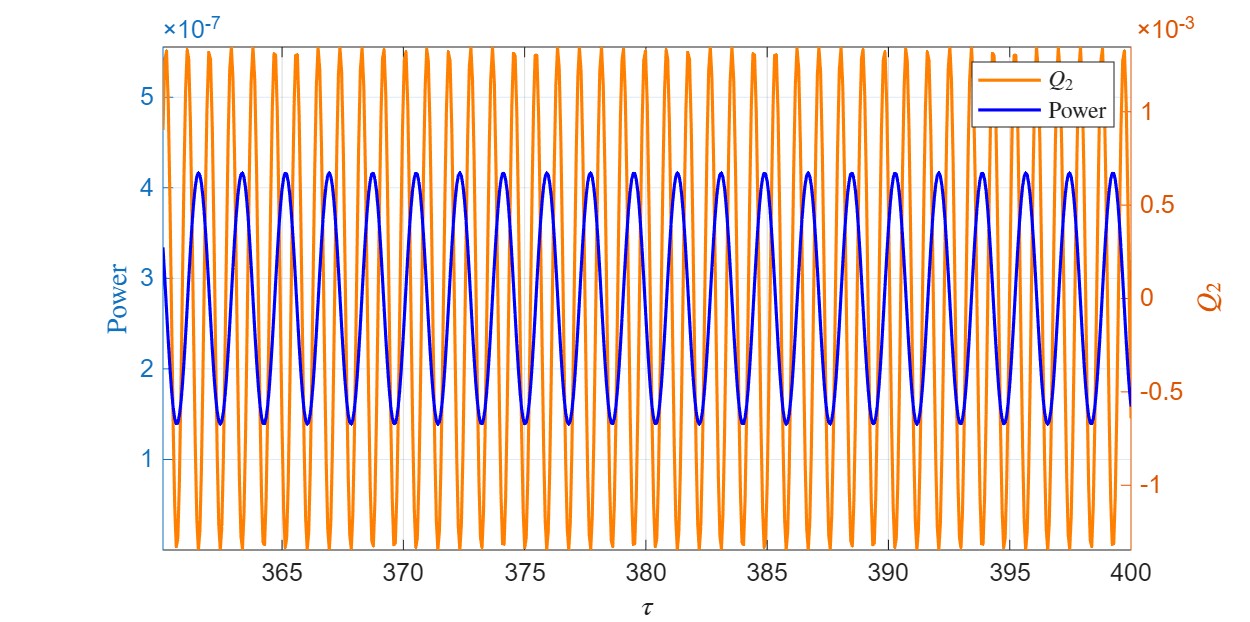}}
\caption{Power harvesting under internal resonance in the new system equation as a function of $\tau$}
\label{fig7}
\end{figure}
The equilibrium points of the equations satisfy $p_1' = p_2' = p_3' = \gamma_1' = \gamma_2' = \gamma_3' = 0$. From the system of equations (\ref{12}), it follows that this nonlinear magnetic levitation system will satisfy:
\begin{equation}
\begin{aligned}
&E^2 = \alpha_5^2 p_2^2 W_4 + \alpha_4^2 p_1^2 + 4\sigma_1^2 W_4 p_2^2 + \frac{2\alpha_4}{\alpha_1 \beta_2} \left[ \right.  \alpha_5 \alpha_2 \beta_1 \beta_3 p_3^3 \\& + \alpha_5 \alpha_3 \beta_2 p_1^2 + 4\sigma_1 \alpha_2 p_3^2 \beta_3 (\sigma_1 + \sigma_2 - \sigma_3)- 4\sigma_1 \beta_2 p_1^2 ( \sigma_1 + \sigma_2 - \frac{3}{8} W_3 p_1 \overline{p_1} ) \left. \right]
\end{aligned}
\end{equation}
where:
\begin{equation}
p_2 = \frac{\sqrt{\left( \alpha_2 \beta_1 \beta_3 p_3^2 + \alpha_3 p_1^2 \beta_2 \right)^2 + 4 \left( \alpha_2 p_3^2 \beta_3 (\sigma_1 + \sigma_2 - \sigma_3) - \ \beta_2 p_1^2 \left( \sigma_1 + \sigma_2 - \frac{3}{8} W_3 p_1 \overline{p_1} \right) \right)^2}}{\beta_2 p_1 \alpha_1 \sqrt{W_4}}
\end{equation}

\begin{equation}
p_3 = \frac{\beta_2 p_1}{\sqrt{\beta_1^2 \beta_3 + 4\beta_3 (\sigma_1 + \sigma_2 - \sigma_3)^2}}
\end{equation}

To investigate the stability of the equations, we assume:
$$
p_i = p_{i0} + p_{i1}, \quad \gamma_i = \gamma_{i0} + \gamma_{i1} \quad (i=1,2,3)
$$
\begin{figure}[t]
\centering
\centerline{\includegraphics[width=1.0\linewidth, height=4cm]{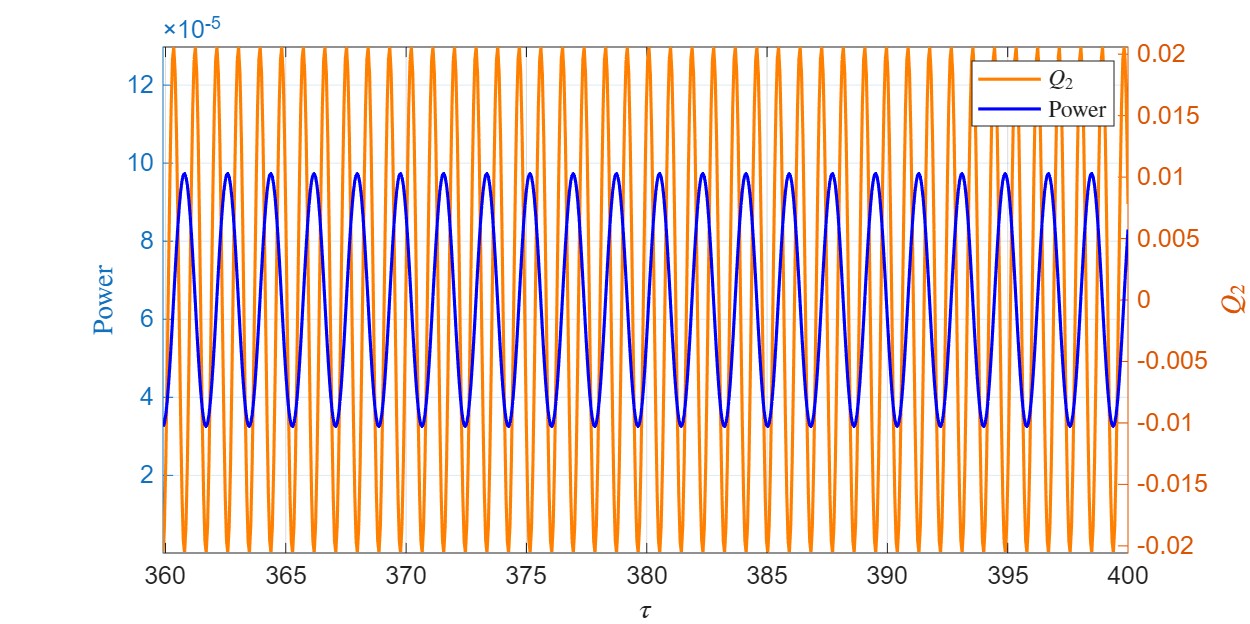}}
\caption{Power harvesting under primary resonance in the new system equation as a function of $\tau$}
\label{fig8}
\end{figure}
Where $p_{i0}$ and $\gamma_{i0}$ represent the equilibrium points of the individual equations in equation (\ref{12}), we can obtain the Jacobian matrix of the system of equations satisfying:
\begin{equation*}
\begin{bmatrix}
p_1' \\
\gamma_1' \\
p_2' \\
\gamma_2' \\
p_3' \\
\gamma_3'
\end{bmatrix}
= J \begin{bmatrix}
p_1 \\
\gamma_1 \\
p_2 \\
\gamma_2 \\
p_3 \\
\gamma_3
\end{bmatrix}
= \begin{bmatrix}
\eta_{11} & \eta_{12} & \eta_{13} & \eta_{14} & \eta_{15} & \eta_{16} \\
\eta_{21} & \eta_{22} & \eta_{23} & \eta_{24} & \eta_{25} & \eta_{26} \\
\eta_{31} & \eta_{32} & \eta_{33} & \eta_{34} & \eta_{35} & \eta_{36} \\
\eta_{41} & \eta_{42} & \eta_{43} & \eta_{44} & \eta_{45} & \eta_{46} \\
\eta_{51} & \eta_{52} & \eta_{53} & \eta_{54} & \eta_{55} & \eta_{56} \\
\eta_{61} & \eta_{62} & \eta_{63} & \eta_{64} & \eta_{65} & \eta_{66}
\end{bmatrix}
\begin{bmatrix}
p_1 \\
\gamma_1 \\
p_2 \\
\gamma_2 \\
p_3 \\
\gamma_3
\end{bmatrix}
\end{equation*}
A detailed analysis of the Jacobian matrix expansion is provided in Appendix B. If the real parts of all eigenvalues of the given Jacobian matrix are negative, the equilibrium points should be globally stable. Otherwise, they become unstable.

\subsection{Primary Resonance Conditions}

Similar to the internal resonance case, this section continues to employ multiscale analysis to solve differential equations. Here, another set of scaling parameters is introduced. Unlike the previous section, the scaling parameters below are accurate to the second order of $\varepsilon$ at most.
\begin{figure}[t]
\centering
\centerline{\includegraphics[width=1.2\linewidth, height=4cm]{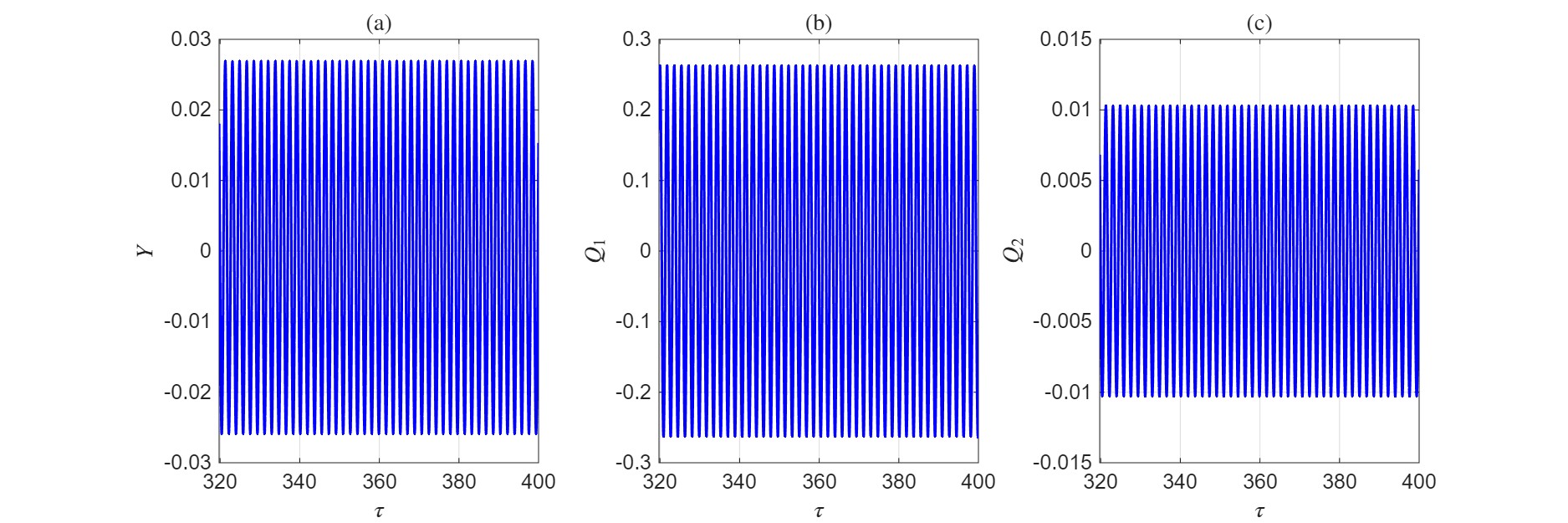}}
\caption{Response of three primary parameters to normalized time $\tau$ under primary resonance conditions}
\label{fig9}
\end{figure}
\begin{equation*}
W_2 = \varepsilon \hat{W_2}, W_3 = \varepsilon \hat{W_3}, \alpha_i = \varepsilon\hat{\alpha_i} \ (i=1,2,3,4,5),E = \varepsilon \hat{E}, \beta_i=\varepsilon\hat{\beta_i} \ (i=1,2)
\end{equation*}
Let: $T_0 = \tau,\quad T_1 = \varepsilon \tau,\quad T_2 = \varepsilon^2 \tau$,and expand $Y,Q_1,Q_2$ as power series in $\varepsilon$:
\begin{equation}
\begin{aligned}
Y &= u_0 (T_0, T_1, T_2) + \varepsilon v_0 (T_0, T_1, T_2) + \varepsilon^2 w_0 (T_0, T_1, T_2) + O(\varepsilon^3)
\\
Q_1 &= u_1 (T_0, T_1, T_2) + \varepsilon v_1 (T_0, T_1, T_2) + \varepsilon^2 w_1 (T_0, T_1, T_2) + O(\varepsilon^3) 
\\
Q_2 &= u_2 (T_0, T_1, T_2) + \varepsilon v_2 (T_0, T_1, T_2) + \varepsilon^2 w_2 (T_0, T_1, T_2) + O(\varepsilon^3) \label{20}
\end{aligned}
\end {equation}
Expand the first- and second-order derivative operators according to the power of $\varepsilon$:
\begin{equation}
\begin{aligned}
\frac{d}{d\tau} &= D_0 + \varepsilon D_1 + \varepsilon^2 D_2 + O (\varepsilon^3)
\\
\frac{d^2}{d\tau^2} &= D_0^2 + 2 \varepsilon D_0 D_1 + \varepsilon^2 (D_1^2 + 2D_0 D_2) + O(\varepsilon^3) \label{21}
\end{aligned}
\end{equation}
$O(\varepsilon^0)$:

\begin{equation}
\begin{aligned}
D_0^2 u_0 + u_0 &= 0 
\\
D_0^2 u_1 + W_4 u_1 &= 0 
\\
D_0^2 u_2 + \beta_3 u_2 &= 0 
\end{aligned}
\end{equation}
\begin{figure}[t]
\centering
\centerline{\includegraphics[width=1.2\linewidth, height=4cm]{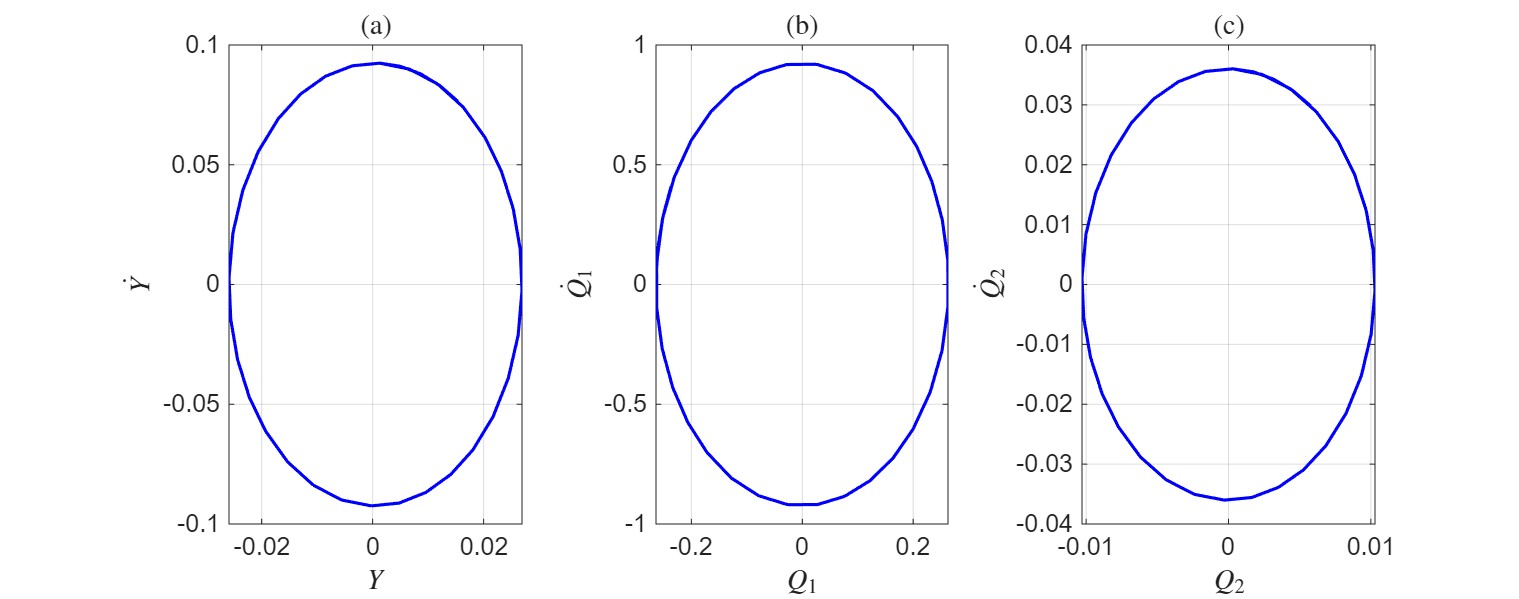}}
\caption{Stable phase plane of three primary parameters under primary resonance conditions}
\label{fig10}
\end{figure}
The final solution yields:
\begin{equation}
\begin{aligned}
u_0 &= A_1 (T_1, T_2) e^{iT_0} + \text{c.c.} 
\\
u_1 &= A_2 (T_1, T_2) e^{i\sqrt{W_4} T_0} + \text{c.c.} 
\\
u_2 &= {A_3(T_1, T_2)} e^{i\sqrt{\beta_3}T_0} + \text{c.c.}  \label{22}
\end{aligned}
\end{equation}
$O(\varepsilon^1)$:
\begin{equation}
\begin {aligned}
(D_0^2 + 1)v_0 &= -2D_0 D_1 u_0 - \alpha_1 D_0 u_1 - \alpha_2 D_0 u_2 \\
(D_0^2 + W_4)v_1 &= -2D_0 D_1 u_1 + \alpha_4 D_0 u_0 \\
(D_0^2 + \beta_3)v_2 &= -2D_0 D_1 u_2 + \beta_2 D_0 u_0\\  \label{24}
\end{aligned}
\end{equation}
Substituting equation (\ref{22}) into equation (\ref{24}), the system of differential equations is ultimately simplified to:

\begin{equation}
\begin{aligned}
(D_0^2 + 1) v_0 &= -2i \frac{\partial A_1}{\partial T_1} e^{iT_0} - i\alpha_1 \sqrt{W_4} A_2 e^{i\sqrt{W_4} T_0} - i\alpha_2 \sqrt{\beta_3} A_3 e^{i\sqrt{\beta_3} T_0} + c.c.\\
(D_0^2 + W_4) v_1 &= -2i \sqrt{W_4} \frac{\partial A_2}{\partial T_1} e^{i\sqrt{W_4} T_0} + i\alpha_4 A_1 e^{iT_0} + c.c.\\
(D_0^2 + \beta_3) v_2 &= -2i \sqrt{\beta_3} \frac{\partial A_3}{\partial T_1} e^{i\sqrt{\beta_3} T_0} + i\beta_2 A_1 e^{iT_0} + c.c.\label{27}
\end{aligned}
\end{equation}
Similar to the internal resonance case, to exclude unbounded solutions and obtain stable solutions, we require: $\frac{\partial A_i}{\partial T_1} = 0 \quad (i=1, 2, 3) \quad \text{i.e.,} \quad A_i = A_i(T_2) $. Solving the system of equations in (\ref{27}), we ultimately obtain:
\begin{figure}[t]
\centering
\centerline{\includegraphics[width=1.2\linewidth, height=5cm]{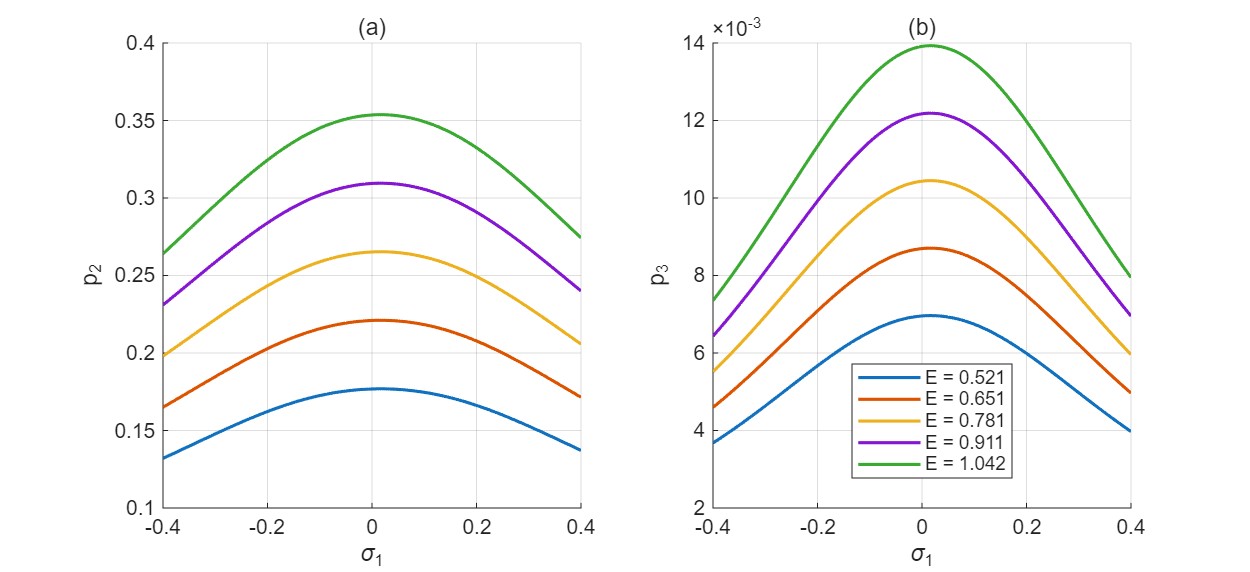}}
\caption{Amplitude of the two primary parameters for different $m$ values under primary resonance conditions as a function of $\sigma_1$}
\label{fig11}
\end{figure}
\begin{equation}
\begin{aligned}
v_0 &= \frac{i\alpha_1 \sqrt{W_4} A_2}{W_4 - 1} e^{i\sqrt{W_4} T_0} + \frac{i\alpha_2 \sqrt{\beta_3} A_3}{\beta_3 - 1} e^{i\sqrt{\beta_3} T_0} + c.c.
\\
v_1 &= \frac{i\alpha_4 A_1}{W_4 - 1} e^{iT_0} + c.c.
\\
v_2 &= \frac{i\beta_2 A_1}{\beta_3 - 1} e^{iT_0} + c.c. \label{29}
\end{aligned}
\end{equation}
$O(\varepsilon^2)$:
\begin{equation}
\begin{aligned}
(D_0^2 + 1) w_0 &= - (2D_0 D_2 + D_1^2) u_0 - 2D_0 D_1 v_0 + W_2 u_0^2 - W_3 u_0^3 \\
&\quad - \alpha_1 D_0 v_1 - \alpha_1 D_1 u_1 - \alpha_2 D_0 v_2 - \alpha_2 D_1 u_2 - \alpha_3 D_0 u_0\\
(D_0^2 + W_4) w_1 &= - (2D_0 D_2 + D_1^2) u_1 - 2D_0 D_1 v_1 \\
&\quad + \alpha_4 D_0 v_0 + \alpha_4 D_1 u_0 - \alpha_5 D_0 u_1 + E \cos\left( {\Omega} T_0 \right)\\
(D_0^2 + \beta_3) w_2 &= - (2D_0 D_2 + D_1^2) u_2 - 2D_0 D_1 v_2 - \beta_1 D_0 u_2 \\
&\quad + \beta_2 D_0 v_0 + \beta_2 D_1 u_0\label{30}
\end{aligned}
\end{equation}
\begin{figure}[t]
\centering
\centerline{\includegraphics[width=1.2\linewidth, height=4cm]{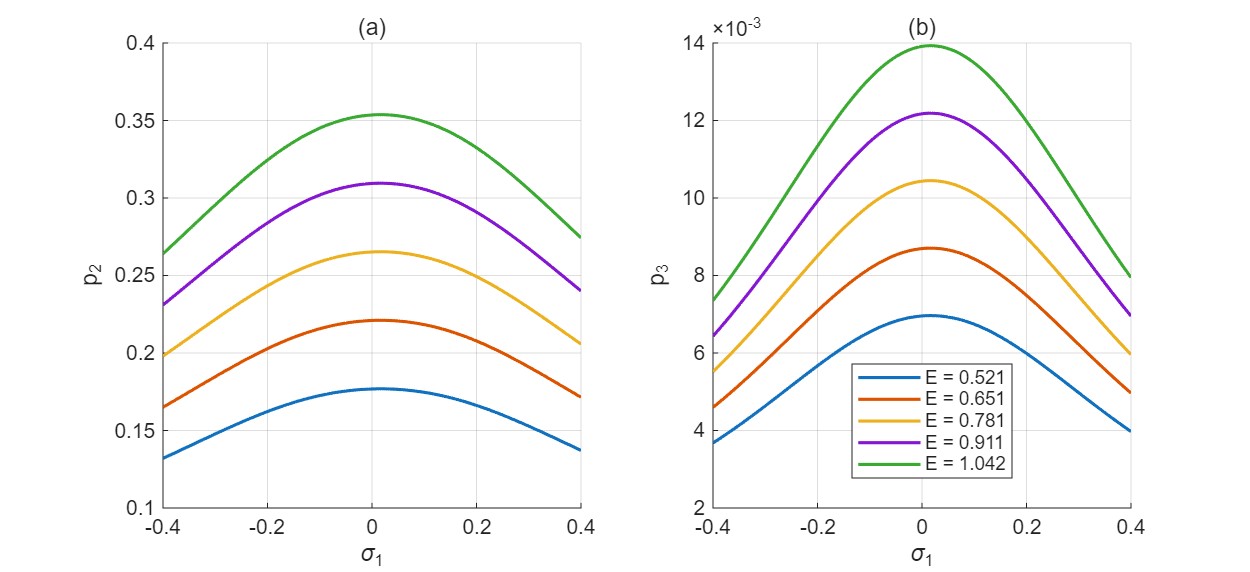}}
\caption{Amplitude of the two primary parameters as $\sigma_1$ varies for different $E$ values under primary resonance conditions}
\label{fig12}
\end{figure}
Substituting equations (\ref{22}, \ref{29}) into the system (\ref{30}) yields:
\begin{equation}
\begin{aligned}
(D_0^2 + 1) w_0 &= -2i A_1' e^{iT_0} + \frac{\alpha_1 \alpha_4 A_1}{W_4 - 1} e^{iT_0} + \frac{\alpha_2 \beta_2 A_1}{\beta_3 - 1} e^{iT_0} - i\alpha_3 A_1 e^{iT_0} \\
&\quad + W_2 A_1^2 e^{2iT_0} + 2W_2 A_1 \overline{A_1} - W_3 A_1^3 e^{3iT_0} - 3W_3 A_1^2 \overline{A_1} e^{iT_0} + c.c.\\
(D_0^2 + W_4) w_1 &= -2i\sqrt{W_4} A_2' e^{i\sqrt{W_4} T_0} - \frac{\alpha_1 \alpha_4 W_4 A_2}{W_4 - 1} e^{i\sqrt{W_4} T_0} - \frac{\alpha_2 \alpha_4 \beta_3 A_3}{\beta_3 - 1} e^{i\sqrt{\beta_3} T_0} \\&- i\alpha_5 \sqrt{W_4} A_2 e^{i\sqrt{W_4} T_0} + \frac{1}{2} E e^{i\widehat{\Omega} T_0} + c.c.\\
(D_0^2 + \beta_3) w_2 &= -2i\sqrt{\beta_3} A_3' e^{i\sqrt{\beta_3} T_0} - i\beta_1 \sqrt{\beta_3} A_3 e^{i\sqrt{\beta_3} T_0} \\& - \frac{\beta_2 \alpha_1 W_4 A_2}{W_4 - 1} e^{i\sqrt{W_4} T_0} - \frac{\beta_2 \alpha_2 \beta_3 A_3}{\beta_3 - 1} e^{i\sqrt{\beta_3} T_0} + c.c.
\end{aligned}
\end{equation}

Similar to the internal resonance case, complex numbers are used here:
\begin{equation}
A_i = \frac{1}{2} p_i e^{iq_i} \quad (i=1,2,3)  \label{35}
\end{equation}

Additionally, two parameters $\sigma_1$ and $\sigma_4$ are introduced, where $\sigma_1$ represents the primary resonance detuning parameter, and $\sigma_4$ characterizes the proximity between the collector circuit frequency and the excitation circuit frequency:
\begin{equation}
\begin{aligned}
\Omega = \sqrt{W_4} + \varepsilon^2 \sigma_1, \quad \sqrt{W_4} = \sqrt{\beta_3} + \varepsilon^2 \sigma_4 \label{36}
\end{aligned}
\end{equation}
To exclude unbounded solutions and obtain stable solutions, the following conditions are required:
\begin{equation}
\begin{aligned}
-&2i A_1' + \frac{\alpha_1 \alpha_4 A_1}{W_4 - 1} + \frac{\alpha_2 \beta_2 A_1}{\beta_3 - 1} - i\alpha_3 A_1 - 3W_3 A_1^2 \overline{A_1} = 0 \\
-&2i\sqrt{W_4} A_2' - \frac{\alpha_1 \alpha_4 W_4 A_2}{W_4 - 1} - \frac{\alpha_2 \alpha_4 \beta_3 A_3}{\beta_3 - 1} e^{i(\sqrt{\beta_3} - \sqrt{W_4}) T_0} - i\alpha_5 \sqrt{W_4} A_2 + \frac{1}{2} E e^{i\sigma_1 T_2} = 0 \\
-&2i\sqrt{\beta_3} A_3' - i\beta_1 \sqrt{\beta_3} A_3 - \frac{\alpha_1 \beta_2 W_4 A_2}{W_4 - 1} e^{i(\sqrt{W_4} - \sqrt{\beta_3}) T_0} - \frac{\alpha_2 \beta_2 \beta_3 A_3}{\beta_3 - 1} = 0 \label{42}
\end{aligned}
\end{equation}
For convenience in subsequent derivations, we introduce two parameters: $\gamma_1 = \sigma_4 T_2 + q_2 - q_3$, $\gamma_2 = \sigma_1 T_1 - q_2$. Substituting equations (\ref{35}, \ref{36}) into (\ref{42}) yields the following system of differential equations:
\begin{equation}
\begin{aligned}
2p_1' &= -\alpha_3 p_1 \\
2p_1 q_1' &= -\left( \frac{\alpha_1 \alpha_4}{W_4 - 1} + \frac{\alpha_2 \beta_2}{\beta_3 - 1} \right) p_1 + \frac{3}{4} W_3 p_1^2 \overline{p_1}\\
2\sqrt{W_4} p_2' &= \frac{\alpha_2 \alpha_4 \beta_3 p_3}{\beta_3 - 1} \sin\gamma_1 - \alpha_5 \sqrt{W_4} p_2 + E \sin\gamma_2 \\
2\sqrt{W_4} p_2 \gamma_2' &= -\frac{\alpha_2 \alpha_4 \beta_3 p_3}{\beta_3 - 1} \cos\gamma_1 + 2\sqrt{W_4} p_2 \sigma_1 - \frac{\alpha_1 \alpha_4 W_4 p_2}{W_4 - 1} + E \cos\gamma_2 \\
2\sqrt{\beta_3} p_3' &= -\beta_1 \sqrt{\beta_3} p_3 - \frac{\alpha_1 \beta_2 W_4 p_2}{W_4 - 1} \sin\gamma_1 \\
2\sqrt{\beta_3} p_3 (\gamma_1' + \gamma_2') &= 2\sqrt{\beta_3} (\sigma_1 + \sigma_4) p_3 - \frac{\alpha_2 \beta_2 \beta_3 p_3}{\beta_3 - 1} - \frac{\alpha_1 \beta_2 W_4 p_2}{W_4 - 1} \cos\gamma_1 \label{45}
\end{aligned}
\end{equation}
The equilibrium points of the equations satisfy $p_2' = \gamma_2' = p_3' = \gamma_1' = 0$ and $p_1=0$. From the system of equations (\ref{45}), it follows that this nonlinear system will satisfy:
\begin{equation}
\begin{aligned}
E =& \sqrt{A^2 + B^2}
\\
A =& \alpha_5 \sqrt{W_4} p_2 + \frac{\alpha_2 \alpha_4 \beta_3 \sqrt{\beta_3} (W_4 - 1) \beta_1 p_3^2}{(\beta_3 - 1) \alpha_1 \beta_2 W_4 p_2}
\\
\quad B = &2\sqrt{W_4} \sigma_1 p_2 - \frac{\alpha_1 \alpha_4 W_4 p_2}{W_4 - 1} - \\
&\frac{(W_4 - 1) \alpha_2 \alpha_4 \beta_3 \left( 2\sqrt{\beta_3} (\sigma_1 + \sigma_4) (\beta_3 - 1) - \alpha_2 \beta_2 \beta_3 \right) p_3^2}{\alpha_1 \beta_2 W_4 p_2(\beta_3-1)^2} \label{40}
\end{aligned}
\end{equation}
and the following equality holds:
\begin{equation}
\frac{\alpha_1 \beta_2 W_4 p_2}{|W_4 - 1|} = \sqrt{\beta_1^2 \beta_3 + \left( 2\sqrt{\beta_3} (\sigma_1 + \sigma_4) - \frac{\alpha_2 \beta_2 \beta_3}{\beta_3 - 1} \right)^2} p_3 \label{41}
\end{equation}
To investigate the stability of the equations, we assume: $p_i = p_{i0} + p_{i1}, \quad \gamma_i = \gamma_{i0} + \gamma_{i1} \quad (i=1,2,3)$,
where $p_{i0}, \gamma_{i0}$ denote the equilibrium points of the respective equations in (\ref{40}). We obtain the Jacobian matrix $J$ of the system of equations satisfying:
\[
\begin{bmatrix}
p_2' \\
\gamma_2' \\
p_3' \\
\gamma_1'
\end{bmatrix}
= J \begin{bmatrix}
p_2 \\
\gamma_2 \\
p_3 \\
\gamma_1
\end{bmatrix}
= \begin{bmatrix}
\delta_{11} & \delta_{12} & \delta_{13} & \delta_{14} \\
\delta_{21} & \delta_{22} & \delta_{23} & \delta_{24} \\
\delta_{31} & \delta_{32} & \delta_{33} & \delta_{34} \\
\delta_{41} & \delta_{42} & \delta_{43} & \delta_{44}
\end{bmatrix}
\begin{bmatrix}
p_2 \\
\gamma_2 \\
p_3 \\
\gamma_1
\end{bmatrix}
\]
Detailed parameter data for the Jacobian matrix can be found in Appendix B.

\section{Results and Discussion}

Based on the model parameters determined in the theoretical study of \cite{16} as cited in \cite{11}, this paper employs the following model parameters: \( W_2 = 1.5000 \), \( W_3 = 1.0000 \), \( W_4 = 2.2783 \), \( \alpha_1 = 0.3247 \), \( \alpha_2 = 0.3247 \), \( \alpha_3 = 0.3125 \), \( \alpha_4 = 0.3248 \), \( \alpha_5 = 0.8400 \), \( \beta_1 = 0.8333 \), \( \beta_2 = 0.3248 \), \( \beta_3 = 2.2783 \), \( E = 0.7812 \), \( \Omega = 3.5 \). This section first elaborates on the improvements and advancements of the new model through phase diagrams, Poincare maps, bifurcation diagrams, and related content; Subsequently, based on the model parameter data, numerical simulations are employed to conduct an in-depth investigation of the magnetic levitation model \textbf{by varying parameters such as mass \(m\), voltage source amplitude \(E\), and damping coefficient \(C_{me}\), the effects of different parameters on system performance are examined}, particularly the stability behavior and energy harvesting efficiency under both internal resonance and primary resonance conditions.

\subsection{Model Enhancement and Advancement}
This paper modifies the energy harvesting circuit of the original model from a purely resistive circuit governed by Kirchhoff's laws to an RLC circuit. \textbf{Compared to the original model, this modification provides a more refined characterization of energy harvesting behavior and better aligns with practical conditions}.
\begin{figure}[t]
\centering
\centerline{\includegraphics[width=1.0\linewidth, height=12cm]{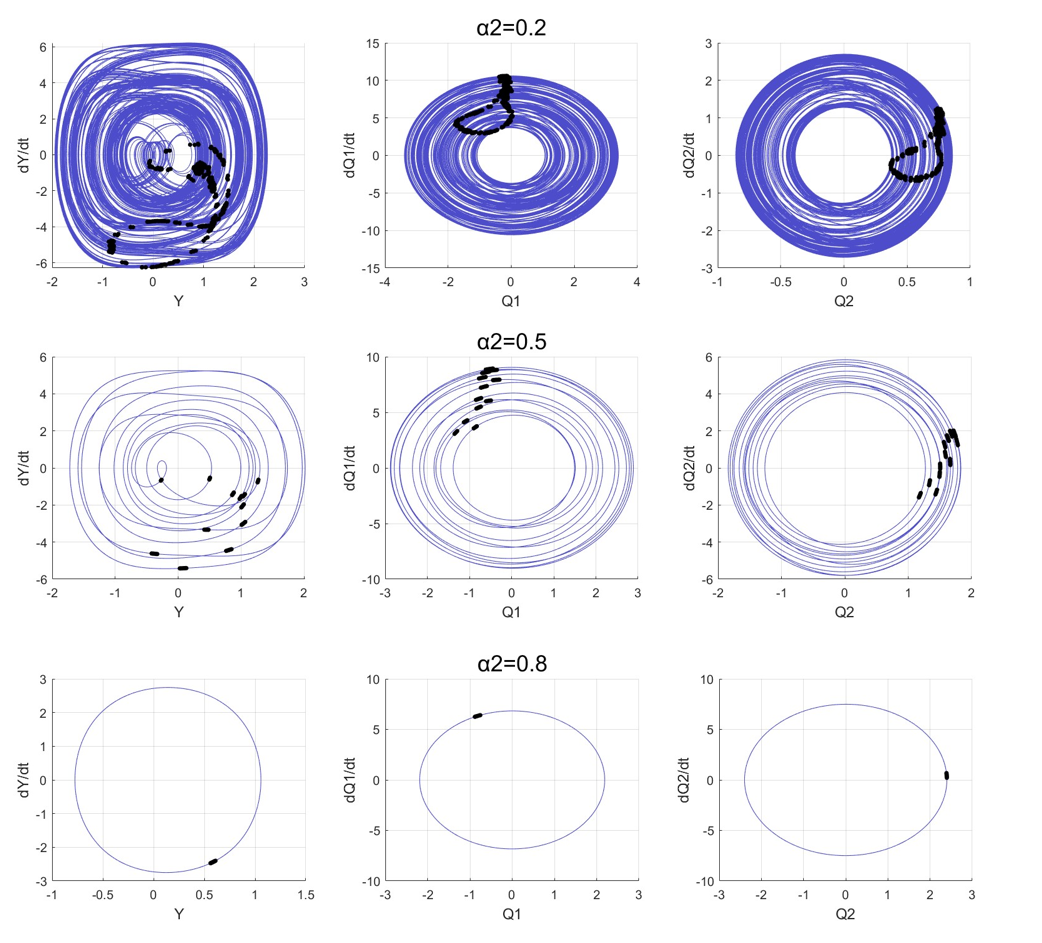}}
\caption{Phase diagrams of $Y, Q_1, Q_2$ (blue curves) and their Poincare maps (black scatter points) for $\alpha_2$ values of $0.2, 0.5, 0.8$}
\label{fig14}
\end{figure}

Regarding energy harvesting capability, Figure \ref{fig文献功率} illustrates the power output from the harvesting circuit as a function of normalized time $\tau$ for both the old system equation's internal resonance and the new system equation's primary resonance. The parameters used correspond to the data employed in this paper. The average collected power calculated for the old system equations is $7.15 \times 10^{-5}$, while the new system equations in this paper yield an average collected power of $6.52 \times 10^{-5}$ using the same experimental data (here, average power is calculated as the integral of steady-state power over time divided by time). Comparing the two curves in Figure \ref{fig文献功率}, although the collection power of the new system equation shows a slight decrease relative to the original system, \textbf{this is due to the further refinement of the energy harvesting circuit in this paper}. The coupling between the energy harvesting circuit and the magnets reduces the amplitude of the intermediate magnet's motion, leading to a slight decrease in collection power. This does not indicate insufficient energy harvesting efficiency in the proposed system but rather an inevitable consequence of system refinement.\textbf{The research value lies in demonstrating that, in practical applications, controlling the capacitance $C_t$ of the harvesting circuit enables frequency matching between the circuit and the mechanical system's vibration frequency ($\beta_2\dot{Y}$), thereby maximizing collected power}. In fact, assuming the mechanical vibration follows a sinusoidal curve with constant amplitude and disregarding parasitic losses from capacitance and resistance, when the collection circuit frequency perfectly matches the mechanical system's vibration frequency, the collected current amplitude depends solely on the collection circuit resistance. It equals the current amplitude in a purely resistive circuit. That is:
\begin{equation*}
\begin{aligned}
   & L\ddot{q}+R\dot{q}+\frac{q}{C}=A\sin \left( \omega _0t \right) \ \left( \omega _{0}^{2}=LC \right) \\
    &\dot{q}=\frac{dq_h\left( t \right)}{dt}+\frac{A}{R}\cos \left( w_0t \right) 
\end{aligned}
\end{equation*}
Where $\frac{dq_h\left( t \right)}{dt}$ arises from electromagnetic damping in the RLC circuit and represents the general solution of the differential equation. \textbf{Note that energy harvesting efficiency peaks during mechanical vibrations of identical period, which inspires the stability analysis of system energy harvesting in this paper.} Regarding this section, this paper characterizes the response of the harvesting circuit to the system by \textbf{introducing detuning parameters $\sigma_3$ and $\sigma_4$, and employs multiscale analysis methods for derivation}.
\begin{figure}[t]
\centering
\centerline{\includegraphics[width=1.2\linewidth, height=8cm]{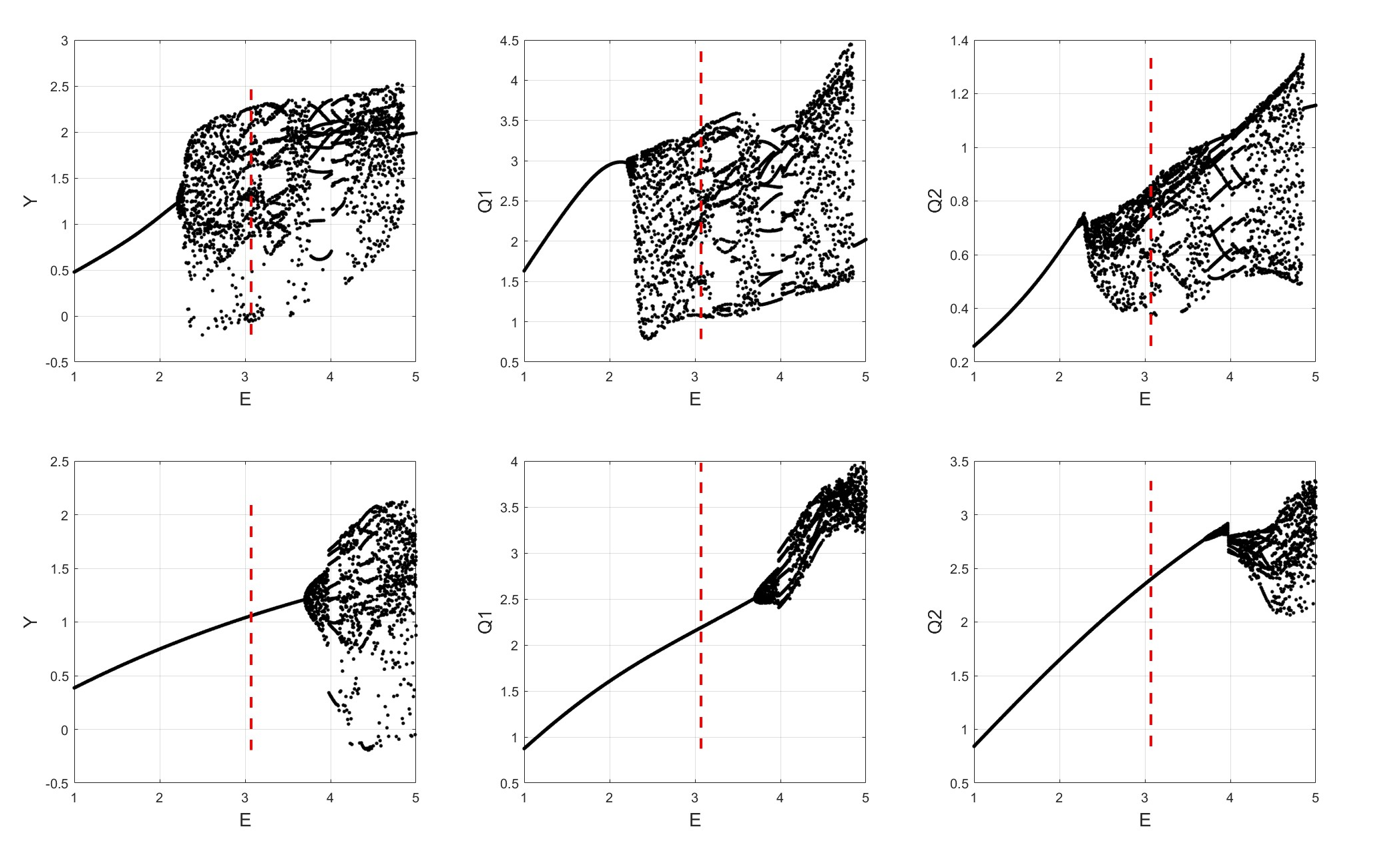}}
\caption{Bifurcation diagram of $Y, Q_1, Q_2$ with respect to $E$ (top figure for $\alpha_2=0.2$, bottom figure for $\alpha_2=0.8$)}
\label{fig16}
\end{figure}
Regarding energy harvesting stability, this paper mitigates the impact of chaotic phenomena in certain voltage ranges on energy harvesting behavior through model refinement. Parameters from reference \cite{17} were substituted into the new model, with test parameters as follows: \( W_2 = 2.0000 \), \( W_3 = 4.0000 \), \( W_4 = 9.108 \), \( \alpha_1 = 0.64944 \), \( \alpha_3 = 0.61996 \), \( \alpha_4 = 0.3248 \), \( \alpha_5 = 0.1499 \), \( \beta_1 = 0.1000 \), \( \beta_3 = 9.108 \), \( E = 3.07 \), \( \Omega = 3.1215 \). For parameters \(\alpha_2\) and \(\beta_2\), take values: $0.0000, 0.2000, 0.5000, 0.8000$ and $0.0000, 0.1500, 0.3750, 0.6000$. When both $\alpha_2$ and $\beta_2$ are 0, the system exhibits the multi-periodic or chaotic behavior of magnet displacement observed in [17]. \textbf{Observing Figure \ref {fig14}. As $\alpha_2$ and $\beta_2$ gradually increase, the phase diagrams of $Y, Q_1, Q_2$ evolve toward limit cycles, and the set of points in the Poincare map gradually concentrates from a dispersed, chaotic pattern}. This indicates that the system's original chaotic behavior is mitigated or even eliminated to some extent. The key parameter influencing the system's energy harvesting stability is $\alpha_2$, the dimensionless coupling coefficient $S_2$ between the energy harvesting circuit and the mechanical circuit. It can be observed that as the coupling coefficient increases, the efficiency of mechanical-to-electrical energy conversion between the mechanical circuit and the energy harvesting circuit improves. \textbf{The periodic harvesting behavior of the harvesting circuit feeds back to the mechanical circuit, causing it to exhibit stronger periodicity, thereby mitigating the negative impact of chaotic behavior on the system.} Observing Figure \ref{fig16}, it is found that when $\alpha_2$ increases from $0.2$ to $0.8$,$E$'s bifurcation critical points with respect to $Y, Q_1, Q_2$ all shift to the right overall. \textbf{This indicates that at higher voltages $E \geq 2.5$, appropriately increasing the coupling coefficient $S_2$ helps improve the system's energy harvesting stability}.

\subsection{Internal Resonance Conditions}
Adjust the capacitance $C_s$ in the excitation circuit to $2.278C_s$ based on the set values of physical parameters, so that the natural frequencies of Equation (\ref{2}) and Equation (\ref{1}) are approximately equal. and adjust the capacitance $C_t$ in the collection circuit to $2.278C_t$ to approximate the natural frequencies of equations (\ref{3}) and (\ref{1}), thereby achieving an internal resonance state. Figure \ref{fig2} displays the periodic time-domain responses of the intermediate magnet displacement, excitation circuit charge, and collection circuit charge during internal resonance. Figure \ref{fig3} displays the stable phase plane of the intermediate magnet displacement, excitation circuit charge, and collection circuit charge during internal resonance. Figure \ref{fig7} illustrates the power output by the collection circuit as a function of normalized time $\tau$ during internal resonance. It can be observed that the magnet exhibits stable periodic motion, and the harvesting circuit outputs stable power. Since the objective of this paper is to enhance energy harvesting efficiency, this section will primarily analyze the charge quantity of the harvesting circuit.

\subsubsection{Analysis of Parameter Effects}

Figure \ref{fig4} illustrates how the amplitudes of three primary parameters vary with $\sigma_1$ for different masses ($0.06kg \leq m \leq 0.14kg$) under internal resonance conditions. Figure \ref{fig4}a represents the displacement of the intermediate magnet, Figure \ref{fig4}b represents the amplitude of the charge in the excitation circuit, and Figure \ref{fig4}c represents the amplitude of the charge in the collection circuit. It can be observed that all amplitude parameters reach their maximum values near $\sigma_1=0$ and exhibit a decreasing trend as the absolute value of $\sigma_1$ increases ($-0.4 \leq \sigma_1 \leq 0.4$). The excitation circuit charge amplitude exhibits a bimodal phenomenon at larger masses, possibly due to bifurcation occurring twice \cite{18}. At $C_{me}=0.3$, a slight bifurcation also appears in the excitation circuit charge amplitude. As $m$ increases ($0.06kg \leq m \leq 0.14kg$), the peak amplitude of the collected charge increases, with the corresponding $\sigma_1$ value at the maximum also rising. Overall, mass exerts a relatively minor influence on the peak amplitude of collected charge. Even when mass increases from 0.06 to 0.14—a 133\% rise—the peak amplitude increases by only 1.2\%. However, as the absolute value of $\sigma_1$ gradually increases, mass's impact on amplitude becomes progressively more significant. For practical systems, this indicates that energy harvesting power is subject to complex mass effects. Moreover, the primary resonance degree—specifically, the proximity between the parameter $\sqrt{W_4}$ and $\Omega$—significantly influences the extent to which mass affects energy harvesting power.

Figure \ref{fig5} illustrates how the amplitudes of the three primary parameters vary with $\sigma_1$ for different $C_{me}$ values ($0.4 \leq C_{me} \leq 0.7$) under internal resonance conditions. Observing Figure \ref{fig5}c, we find that $C_{me}$ significantly affects the amplitude of the charge quantity in the energy harvesting circuit. As $C_{me}$ decreases ($0.4 \leq C_{me} \leq 0.7$), the peak amplitude of the charge quantity increases. The amplitude peaks near $\sigma_1=0$. As the absolute value of $\sigma_1$ increases gradually ($-0.4 \leq \sigma_1 \leq 0.4$), the amplitude of the charge quantity begins to decrease sharply, and the influence of $C_{me}$ on the amplitude diminishes progressively. For smaller damping coefficients ($C_{me} \leq 0.4$), the amplitude of the central magnet's motion increases sharply, indicating severe oscillations in the magnetic levitation system. This also demonstrates that energy harvesting power is affected by system damping, with the magnitude of this effect peaking during primary resonance ($\sqrt{W_4} \approx \Omega$). When the system is no longer in primary resonance, this influence diminishes sharply.

Figure \ref{fig6} illustrates how the amplitudes of three key parameters vary with $\sigma_1$ for different $E$ values ($0.3 \leq E \leq 1$) under internal resonance conditions. Observing Figure \ref{fig6}c reveals that $E$ significantly impacts the amplitude of the energy harvesting circuit's charge quantity. Adjusting the normalized voltage $E$ from 0.3 to 1 increases the peak amplitude of the energy harvesting circuit charge from 0.084 to 0.252. This indicates that energy harvesting power is influenced by the normalized voltage (i.e., excitation voltage) ($0.3 \leq E \leq 1$), with higher normalized voltages yielding greater harvesting power.

The above analysis of internal resonance shares a common feature: when the system is near primary resonance, i.e., around $\sigma_1=0$, the energy harvesting power of the collection circuit shows a significant increase.

\subsubsection{Analysis of Stability}

Under conditions of high excitation voltage and low mass ($e\geq15V,m\leq0.01kg$), the system exhibits greater instability. For medium mass ($m=0.1kg$), high excitation voltage has a relatively minor impact on system stability. Therefore, appropriately increasing the mass $m$ can simultaneously enhance both energy harvesting efficiency and stability. However, excessively large masses carry the risk of collision between the center magnet and the bottom magnet \cite{18}. The mechanical-electrical damping parameter $C_{me}$ has a negligible effect on system stability.

\subsection{Primary Resonance Conditions}

Based on the previously observed enhancement of power collection by primary resonance, the capacitance $C_s$ in the excitation circuit was adjusted to $0.186C_s$ to approximate the natural frequency and normalized frequency $\Omega$ in Equation (2), and adjusting the capacitance $C_t$ in the collection circuit to $0.186C_t$ makes the natural frequency and normalized frequency $\Omega$ in Equation (3) approximately equal, thereby achieving the primary resonance state ($\sigma_1,\sigma_4\approx0$). Figure \ref{fig7} shows the power output from the harvesting circuit as a function of normalized time $\tau$ during internal resonance. Figure \ref{fig8} shows the power output from the harvesting circuit as a function of normalized time $\tau$ during primary resonance. Comparing Figure \ref{fig7} and Figure \ref{fig8} reveals that \textbf{the average power harvested by the circuit during primary resonance increases by two orders of magnitude compared to internal resonance}. This demonstrates that the novel system achieves a significant enhancement in power harvesting under primary resonance conditions, representing an important complement to current energy harvesting systems that rely on internal resonance. Figure \ref{fig9} displays the periodic time-domain responses of the intermediate magnet displacement, excitation circuit charge, and harvesting circuit charge during primary resonance. Figure \ref{fig10} illustrates the stable phase plane for all three quantities during primary resonance.

\subsubsection{Analysis of Parameter Effects}

Figure \ref{fig11} illustrates how the two amplitudes vary with $\sigma_1$ under primary resonance conditions for different masses ($0.05kg \leq m \leq 0.25kg$). Observing Figure \ref{fig11}b reveals that $m$ significantly influences the amplitude of the energy harvesting circuit's charge quantity. As $m$ decreases ($0.05kg \leq m \leq 0.25kg$), the peak amplitude of the harvesting circuit's charge increases, while the amplitude of the excitation circuit's charge quantity is less affected, with peak differences between different masses being less than 5\%. The amplitude of collected charge reaches its peak near $\sigma_1=0$. As the absolute value of $\sigma_1$ gradually increases ($-0.6 \leq \sigma_1 \leq 0.6$), the amplitude begins to decrease sharply toward zero, and the influence of $m$ on this behavior diminishes progressively. This indicates that energy harvesting power is influenced by mass and primary resonance frequency ($0.05kg \leq m \leq 0.25kg$). A smaller central magnet mass and a closer parameter $\sqrt{W_4}$ to $\Omega$ both significantly enhance the harvested power magnitude, consistent with prior analysis.

Figure \ref{fig12} illustrates how the two amplitudes vary with $\sigma_1$ under primary resonance for different normalized voltages ($0.521 \leq E \leq 1.042$). Adjusting $E$ from 0.521 to 1.042 causes the peak amplitude of the collected circuit charge to vary from 0.007 to 0.014, representing a 100\% increase. This indicates that energy harvesting power is influenced by the normalized voltage (i.e., excitation voltage) within the range $0.521 \leq E \leq 1.042$, with higher voltages yielding greater power collection. Additionally, it is observed that the approximate solution does not include the normalization coefficient $\alpha_3$. This indicates that the mechanical-electrical damping parameter $C_{me}$ exerts a negligible influence on the system during primary resonance, rendering it negligible under the assumption of approximate solution.

\subsubsection{Analysis of Stability}

At low voltages ($0.5V \leq e \leq 6V$), the system remains relatively stable under primary resonance conditions. For medium to high voltages ($e \geq 15V$), the stable region shrinks dramatically, with stability primarily observed in the range $-0.5 \leq \sigma_1 \leq 0.2$ and simultaneously $-0.5 \leq \sigma_2 \leq 0$. The mechanical-electrical damping parameter $C_{me}$ has a relatively minor impact on system stability. The system remains stable when the mass $m$ is in the range $0.06 \leq \sigma_2 \leq 1.0$. However, when the mass is excessively small ($m \leq 0.01kg$), the unstable region actually expands. Therefore, maintaining the mass $m$ at a moderate level is beneficial for system stability.

\section{Conclusions}
This paper systematically investigates the influence of mechanical and electrical parameters on the system dynamics and energy harvesting performance by establishing and analyzing a nonlinear magnetic levitation system model that comprehensively considers the coupling between the excitation circuit, energy harvesting circuit, and mechanical system. \textbf{Using a multiple scales method, approximate analytical solutions were derived for both internal resonance and primary resonance scenarios, with theoretical predictions validated through numerical simulations}. The main conclusions are as follows. \\
Internal Resonance Condition:
\begin{itemize}
\item 
By adjusting capacitance to $2.278C_s$ and $2.278C_t$, the system enters internal resonance, where the central magnet, excitation circuit charge, and harvesting circuit charge all exhibit stable periodic motion.
\item 
Within the range $0.06kg \leq m \leq 0.14kg$, mass has a negligible effect on the peak collected charge, but its influence becomes increasingly significant as the frequency deviates from the primary resonance frequency.
\item 
Within the ranges $0.4 \leq C_{me} \leq 0.7$ and $0.3 \leq E \leq 1$, the damping coefficient and normalized voltage significantly affect the amplitude of collected charge. Their influence peaks near the resonance point and rapidly diminishes as frequency deviates.
\end{itemize}
Primary Resonance Condition:
\begin{itemize}
\item 
By further adjusting the capacitance to $0.186C_s$ and $0.186C_t$ to satisfy the primary resonance condition, the system exhibits more stable periodic motion with reduced fluctuations in collected power.
\item 
\textbf{Under primary resonance, the average collected power increased by two orders of magnitude compared to internal resonance}, indicating this state significantly enhances energy harvesting efficiency and serves as a complementary approach to internal resonance-based collection.
\item 
Within the ranges $0.05kg \leq m \leq 0.25kg$ and $0.521 \leq E \leq 1.042$, reducing mass and increasing normalized voltage further enhances the collected power amplitude, with the most pronounced effect near the resonance point. The damping coefficient has negligible influence under primary resonance conditions.
\end{itemize}

\newpage
\appendix
\section{Model Supplement}
\label{app1}
The following provides supplementary explanations for parameters not appearing in the paper:
\begin{align*}
C_{me} = C_m + C_e, \
C_e = \frac{S_2^2}{R_t}, \
R_t = R_{\text{int}} + R_{\text{load}}, \
S_2 = nBL
\end{align*}
where $R_{\text{int}}$ denotes the internal resistance of the coil, $R_{\text{load}}$ represents the load resistance of the energy harvesting circuit, $C_m$ refers to the mechanical damping coefficient, and $C_e$ denotes the electrical damping coefficient.
Similar to Reference \cite{16}, kinetic energy is given by the motion of the central magnet and the inductances of the two circuits. Explanations for the variables will be provided uniformly after Equation \ref{A}:
\begin{equation*}
T=\frac{1}{2}(m\dot{x}^2+L_s\dot{q}_1^2+L_t\dot{q}_2^2)
\end{equation*}
Potential energy is given by the gravitational potential energy of the central magnet, the restoring forces of the bottom and top coils, the coupling between the displacement of the central magnet and the two circuits, and the capacitance of the two circuits:
\begin{equation*}
V = mgx + \frac{1}{2}x^2 (K_{\text{magb}} - K_{\text{magt}}) + S_1\dot{q}_1x + S_2\dot{q}_2x+\frac{1}{2} \frac{q_1^2}{C_s}+\frac{1}{2} \frac{q_2^2}{C_t}
\end{equation*}
Similarly, the dissipation function is obtained as:
\begin{equation*}
D = \frac{1}{2} \left[ c_{\text{me}} \dot{x}^2+ R_s \dot{q}_1^2+ R_t \dot{q}_2^2 \right]
\end{equation*}
Given the energy relations of the system, the Lagrangian function defined by the difference between kinetic and potential energy is obtained:
\[
\begin{split}
\begin{aligned}
L =& \frac{1}{2} c_{\text{me}} \dot{x}^2 + \frac{1}{2} L_s \dot{q}_1^2 + \frac{1}{2} L_t \dot{q}_2^2 - mgx - \frac{1}{2} x^2 (K_{\text{magb}}) - K_{\text{magt}}) \\&- S_1\dot{q}_1x - S_2\dot{q}_2x-\frac{1}{2} \frac{q_1^2}{C_s}-\frac{1}{2} \frac{q_2^2}{C_t}
\end{aligned}
\end{split}
\]
Using the Euler-Lagrange equations:
\[
\frac{\mathrm{d}}{\mathrm{d}t} \left( \frac{\partial L}{\partial \dot{q}_i} \right) - \frac{\partial L}{\partial q_i} + \frac{\partial D}{\partial \dot{q}_i} = F_{\text{ext} \, i} \quad (i = x, q_1,q_2)
\]
This ultimately yields equations (\ref{M}–\ref{A}) from the preceding section, where the voltage in equation (\ref{N}) is provided by the term $e\cos(\widehat{\Omega}t)$.
\newpage
\section{Jacobian Matrix}
\label{app2}
\[
J_1=\begin{pmatrix}
\zeta _1&		\zeta _2&		\zeta _3&		\zeta _4&		\zeta _5&	\zeta _6\\
\end{pmatrix} 
\]
\begin{align*}
\zeta_1 = \quad & \left[ -\frac{\alpha_3}{2},\ 
\frac{\sigma_1 + \sigma_2}{p_1},\ 
\frac{\alpha_4 \cos\gamma_1}{2\sqrt{W_4}},\ 
\frac{\alpha_4 \sin\gamma_1}{2p_2 \sqrt{W_4}},\ 
\frac{\beta_2 \cos\gamma_3}{2\sqrt{\beta_3}},\ 
\frac{\sigma_1 + \sigma_2 - \sigma_3}{p_1} \right]^T \\
\zeta_2 = \quad & \left[ \frac{1}{2}\alpha_1 \sqrt{W_4}\, p_2 \sin\gamma_1,\ 
\frac{1}{2p_1}\alpha_1 \sqrt{W_4}\, p_2 \cos\gamma_1,\ 
-\frac{\alpha_4 p_1 \sin\gamma_1}{2\sqrt{W_4}},\ 
\frac{\alpha_4 p_1 \cos\gamma_1}{2p_2 \sqrt{W_4}},\ 
0,\ 
0 \right]^T \\
\zeta_3 = \quad & \left[ -\frac{1}{2}\alpha_1 \sqrt{W_4}\cos\gamma_1,\ 
\frac{1}{2p_1}\alpha_1 \sqrt{W_4}\sin\gamma_1,\ 
-\frac{\alpha_5}{2},\ 
\frac{\sigma_1}{p_2},\ 
0,\ 
0 \right]^T \\
\zeta_4 = \quad & \left[ 0,\ 
0,\ 
\frac{E \cos\gamma_2}{2\sqrt{W_4}},\ 
-\frac{E \sin\gamma_2}{2p_2 \sqrt{W_4}},\ 
0,\ 
0 \right]^T \\
\zeta_5 = \quad & \left[ -\frac{1}{2}\alpha_2 \sin\gamma_3,\ 
-\frac{\alpha_2 \cos\gamma_3}{2p_1},\ 
0,\ 
0,\ 
-\frac{\beta_1}{2},\ 
-\frac{\sigma_1 + \sigma_2 - \sigma_3}{p_3} \right]^T \\
\zeta_6 = \quad & \left[ -\frac{1}{2}\alpha_2 p_3 \cos\gamma_3,\ 
\frac{\alpha_2 p_3 \sin\gamma_3}{2p_1},\ 
0,\ 
0,\ 
p_3 (\sigma_1 + \sigma_2 - \sigma_3),\ 
-\frac{\beta_1}{2} \right]^T
\end{align*}

\vspace {15mm}
\[
J_2=\begin{pmatrix}
\xi _1&		\xi _2&		\xi_3&		\xi _4\\
\end{pmatrix} 
\]
\begin{align*}
\xi_1 = \quad & \left[ -\dfrac{\alpha_5}{2},\ 
\sigma_1 - \dfrac{\alpha_1 \alpha_4 \sqrt{W_4}}{2(W_4 - 1)},\ 
-\dfrac{\alpha_1 \beta_2 W_4 \sin\gamma_1}{2\sqrt{\beta_3}\, (W_4 - 1)},\ 
-\dfrac{\alpha_1 \beta_2 W_4 \cos\gamma_1}{2\sqrt{\beta_3}\, p_3 (W_4 - 1)} \right]^T \\
\xi_2 = \quad & \left[ \dfrac{E \cos\gamma_2}{2\sqrt{W_4}},\ 
-\dfrac{E \sin\gamma_2}{2\sqrt{W_4}\, p_2},\ 
0,\ 
0 \right]^T \\
\xi_3 = \quad & \left[ \dfrac{\alpha_2 \alpha_4 \beta_3 \sin\gamma_1}{2\sqrt{W_4}\, (\beta_3 - 1)},\ 
-\dfrac{\alpha_2 \alpha_4 \beta_3 \cos\gamma_1}{2\sqrt{W_4}\, p_2 (\beta_3 - 1)},\ 
-\dfrac{\beta_1}{2},\ 
\sigma_1 + \sigma_4 - \dfrac{\alpha_2 \beta_2 \sqrt{\beta_3}}{2(\beta_3 - 1)} \right]^T \\
\xi_4 = \quad & \left[ \dfrac{\alpha_2 \alpha_4 \beta_3 p_3 \cos\gamma_1}{2\sqrt{W_4}\, (\beta_3 - 1)},\ 
\dfrac{\alpha_2 \alpha_4 \beta_3 p_3 \sin\gamma_1}{2\sqrt{W_4}\, p_2 (\beta_3 - 1)},\ 
-\dfrac{\alpha_1 \beta_2 W_4 p_2 \cos\gamma_1}{2\sqrt{\beta_3}\, (W_4 - 1)},\ 
\dfrac{\alpha_1 \beta_2 W_4 p_2 \sin\gamma_1}{2\sqrt{\beta_3}\, p_3 (W_4 - 1)} \right]^T
\end{align*}

\end{document}